\documentclass{article}
\usepackage{amsmath}
\usepackage{amsthm}
\usepackage{latexsym}
\usepackage{amssymb}
\usepackage[all]{xy}
\usepackage{graphicx}

\theoremstyle{definition}
\newtheorem{thm}{Theorem}[section]
\newtheorem{defi}[thm]{Definition}
\newtheorem{prop}[thm]{Proposition}
\newtheorem{cor}[thm]{Corollary}

\newtheorem{algo}[thm]{Algorithm}
\newtheorem{lemma}[thm]{Lemma}
\newtheorem{rem}[thm]{Remark}
\newtheorem{example}[thm]{Example}

\newcommand\ZZ{\ensuremath{\mathbb{Z}}}
\newcommand\CC{\ensuremath{\mathbb{C}}}

\newcommand\RR{\ensuremath{\mathbb{R}}}
\newcommand\NN{\ensuremath{\mathbb{N}}}
\newcommand\QQ{\ensuremath{\mathbb{Q}}}

\newcommand\TT{\ensuremath{\mathbb{T}}}
\newcommand\V[1]{\mbox{\boldmath$#1$}}
\newcommand\zet[1]{\left\vert {#1} \right\vert}

\renewcommand\hat{\widehat}
\renewcommand\tilde{\widetilde}

\newcommand\tp[1]{{#1}_{\mathrm{top}}}
\newcommand\kagi[1]{\left\langle {#1}\right\rangle }

\newcommand\im{\mathrm{Im}\,}

\title{Perturbation of Perron roots and The max-plus spectral theorem}
\author{Shinsuke Iwao\footnote{E-mail: iwao@gem.aoyama.ac.jp}\\
College of Science and Engineering, Aoyamagakuin University,\\
5-10-1, Fuchinobe, Cyu\={o}-ku, Sagamihara-shi, Kanagawa Japan 252-5258
}

\begin{document}
\maketitle

{\abstract
In this paper, we consider the Perron theorem over the real Puiseux field.
We introduce a recursive method for calculating Perron roots and Perron vectors of positive Puiseux matrices (which satisfy some condition of genericness) by means of combinatorics based on the tropical linear algebra.
}

\section{Introduction}\label{sec1}

A matrix, and in particular a vector, is said to be {\it positive} ({\it resp.\,nonnegative}) if all its entries are positive ({\it resp.}\ nonnegative).
The first important discovery on the spectral  property of positive matrices was made by O. Perron  \cite{perron1906grundlagen} in 1907.
Perron's results were extended to nonnegative matrices by G. Frobenius \cite{frobenius1909matrizen1,frobenius1909matrizen2,frobenius1912matrizen}, who proved the following:
\begin{thm}[The Perron-Frobenius theorem]\label{frob}
Let $A\in \RR_\geq^{N\times N}$ be an irreducible and nonnegative matrix and $\rho(A)$ be its spectral radius.
Then, there exist a natural number $d\geq 1$ (called the {\it period} of $A$) such that the matrix $A$ satisfies the following properties (i--iv):
(i) $\rho(A)\exp(2\pi\sqrt{-1}\cdot \frac{a}{d})$ is a simple eigenvalue for each $0\leq a<d$.
(ii) A corresponding eigenvector to $\rho(A)$ can be chosen to be entry-wise positive. 
(iii) Any nonnegative eigenvector of $A$ corresponds to $\rho(A)$.
(iv) No eigenvalue $\lambda$, except for $\rho(A)\exp(2\pi\sqrt{-1}\cdot\frac{a}{d})$, satisfies $\zet{\lambda}=\rho(A)$.
$\qed$
\end{thm}
The positive eigenvalue $\rho(A)$ is called the {\it Perron root of $A$}, and its corresponding eigenvector which is entry-wise positive is called a {\it Perron vector of $A$}. 

In this paper, we will be concerned with an analogue of Perron's result for positive matrices over the {\it real Puiseux field} $F:=\bigcup_{p=1}^\infty \RR((t^{1/p}))$.
It is well-known that the analogue of the Perron-Frobenius theorem holds over any {\it real closed field} (especially, over $F$).
(To my best knowledge, the first (non-constructive) proof was given by Eaves, Rothblum and Schneider in 1995 \cite{Eaves1995}\footnote{In \cite{Eaves1995}, they found that the Perron-Frobenius theorem can be expressed by a {\it sentence in the language of ordered fields} (see Appendix \S \ref{secaa}).
From the {\it Tarski-Seidenberg principle} \cite{seidenberg1954new,tarski1948decision}, they concluded that the Perron-Frobenius theorem is {\it provable} in the language of ordered fields.
}).
On the other hand, it is also known that it is a perplexing problem to give an explicit, constructive formulas to compute all the Puiseux series coefficients of $\rho(A)$.
This kind of problem is one of the most important concerns of perturbation theory (see, for example,
the standard texts \cite{baumgartel1985analytic,kato1995perturbation,lancaster1987analytic} of analytic perturbation theory).

%
In the past 10 years, a new tropical-algebraic approach to estimate perturbed eigenvalues and eigenvectors has been developed by the experts in the field of tropical algebra.
It was shown by Akian, Bapat and Gaubertin \cite{akian2004min,Akian2004103} that generically, the valuations of the eigenvalues of a Puiseux matrix (= a matrix whose entries are Puiseux series) coincide with the eigenvalues of the associated tropical matrix.
They also gave a generalization of the classical Lidski\u{i}-Vi\u{s}ik-Ljusternik theorem \cite{lidskii1966perturbation,visik1960}, which provides the leading terms of eigenvalues and eigenvectors of a linearly perturbed matrix $Y=Y_0+tY_1$.
(For other resent developments on this topic, see \cite{akian2014tropical,gaubert2009tropical}.)

Inspired by the work \cite{akian2004min}, we ask whether there exists a recursive method to compute Perron roots and vectors of a generally perturbed matrix
\[
Y=Y_0+tY_1+t^2Y_2+t^3Y_3+\cdots
\]
by means of tropical mathematics.
In this paper, we introduce a combinatorial and constructive proof of the Perron theorem over $F$ under the {\bf condition of genericness} (\S \ref{sec.generic})\footnote{Any linearly perturbed matrix $Y_0+tY_1$ satisfies this condition.}.

For the first example, consider the real Puiseux matrix
$
Y=
\left( 
\begin{array}{cc}
1-t & 1+t \\ 
1+t^2 & 1-t
\end{array} 
\right).
$
It is not hard to find the Perron root $\lambda_0=2$ and its corresponding Perron vector $\V{v}_0={}^T(1,1)$ of the {\it leading term $\tp{Y}=
\left( 
\begin{array}{cc}
1 & 1 \\ 
1 & 1
\end{array} 
\right)
$}.
As $2$ is a simple eigenvalue of $\tp{Y}$, we soon obtain the pair
\[
\lambda=2-\frac{1}{2}t+\frac{3}{8}t^2+\cdots,\qquad 
\V{v}=
\left(
\begin{array}{c}
1 \\ 
1
\end{array} 
\right)
+
\left(
\begin{array}{c}
1/2 \\ 
0
\end{array} 
\right)t
-
\left(
\begin{array}{c}
5/8 \\ 
0
\end{array} 
\right)t^2
+\cdots,
\]
with $Y\V{v}=\lambda\V{v}$ by applying the recursive formula (\ref{eq6}).
Since both $\lambda$ and $\V{v}$ are positive, we conclude that $\lambda$ is the Perron root of $Y$.

When the leading term of the matrix is not {\it irreducible}, for example
$
Y=\left(
\begin{array}{cc}
1-t & 1+t \\ 
t^2 & 1-t
\end{array} 
\right)$
or
$
Y=\left(
\begin{array}{ccc}
1-t & 1+t & t^2\\ 
1+t^2 & 1-t & t\\
t & t^2 & 2t^2
\end{array} 
\right)
$,
the situation is much complicated (see Example \ref{ex2.17}).
The key of our strategy is to find a diagonal matrix $\delta$ such that $\delta^{-1}Y\delta$ has a simpler structure. 

Our method to find such $\delta$ is split into two steps.

By using the {\it max-plus spectral theorem} (Proposition \ref{trop}), which is a tropical analog of the Perron-Frobenius theorem, one can find some diagonal matrix $\delta_1$ such that the matrix $\delta_1^{-1}Y\delta_1$ is a {\it flat-slanted form} (Section \ref{sec4}).
Roughly speaking, under the transformation $Y\mapsto Y':=\delta_1^{-1}Y\delta_1$, the unevenness of the valuations is fixed.

For further calculations, we have to take into account the magnitudes of coefficients.
By using a combinatorial methods on weighted graphs, one can find a diagonal matrix $\delta_2$ such that $\delta_2^{-1}Y'\delta_2$ is a {\it gently-slanted form} (Section \ref{sec5}), that is a positive Puiseux matrix which admits a recursive method for calculating the Perron root.

The main theorem of this article is:
\begin{thm}
Let $Y$ be a square matrix whose entries are positive Puiseux series.
If $Y$ satisfies the {\bf condition of genericness} (\S \ref{sec.generic}), we have the recursive method (Algorithm \ref{algo}) to calculate the Perron root and a Perron vector of $Y$.
Especially, we obtain a constructive proof of the Perron theorem for some class of Puiseux matrices.
\end{thm}


The contents of this paper are as follows:
In Section \ref{sec2}, we introduce a recursive method to calculate the eigenvalues and eigenvectors of a complex Puiseux matrix whose leading term has a semi-simple eigenvalue.
Although there are numerous literatures and sophisticated results on this topic (see, for example \cite{baumgartel1985analytic,hryniv1999perturbation,lidskii1966perturbation,visik1960,welters2011explicit}), we shall use an alternative method which is suitable for dealing with {\it positivity}.
In Sections \ref{sec3}--\ref{sec5}, we discuss about the combinatorial aspects of real Puiseux matrices.
The main aim in these sections is to introduce two canonical forms of real Puiseux matrices, one of which is a {\it flat-slanted form} (Section \ref{sec4})  and the other is a {\it gently-slanted form} (Section \ref{sec5}).
In Section \ref{sec6}, we give a construction of the Perron root/vector of a positive Puiseux matrix.
We also give some examples and applications of our method in Section \ref{sec7}.
%
%

\section{Reviews on spectral problem of complex Puiseux matrices}\label{sec2}

In this section, we briefly review recursive techniques of calculating eigenvalues and eigenvectors of complex Puiseux matrices whose leading term have semi-simple eigenvalues. 
For readers who are interested in general methods in this field, we recommend the textbooks \cite{baumgartel1985analytic,kato1995perturbation,lancaster1987analytic}, and the references therein.

\subsection{Module of approximate eigenvectors}\label{sec2.2}

Let $K:=\bigcup_{p=1}^\infty\CC((t^{1/p}))$ be the complex Puiseux series field and $\TT:=\QQ\cup\{+\infty\}$ be the tropical semifield.
Define the valuation map $\vartheta:K\to \TT$ by $\vartheta(c)=0$ ($c\in\CC^\times$), $\vartheta(t)=1$ and $\vartheta(0)=+\infty$.
Let $R:=\{\V{x}\in K\,\vert\,\vartheta(\V{x})\geq 0 \}$ be the valuation ring of $K$ and $I:=\{\V{x}\in K\,\vert\,\vartheta(\V{x})> 0  \}$ be the unique maximal ideal of $R$.
For $\Lambda\in \QQ$, set 
$$
I[\Lambda]:=\{\V{x}\in K\,\vert\,\vartheta(\V{x})\geq \Lambda \},\qquad
I(\Lambda):=\{\V{x}\in K\,\vert\,\vartheta(\V{x})> \Lambda \}.
$$

Let $Y=(a_{i,j})\in K^{M\times N}$ be a Puiseux matrix.
The {\it valuation of $A$} is a rational number $v(Y):=\min\{\vartheta(a_{i,j})\}$.
As the number of entries are finite, there uniquely exists a rational number $q(Y):=\min\{r\,\vert\, a_{i,j}\in \CC((t^{1/r})),\ \forall i,j\}$.
Therefore, $Y$ admits the expansion
$$
Y=Y_{v}t^{v}+Y_{v+1/q}t^{v+1/q}+Y_{v+2/q}t^{v+2/q}+\cdots,\qquad Y_v\neq O,\ Y_r\in \CC^{M\times N}.
$$
The matrix $\tp{Y}:=Y_v$ is called the {\it leading term} of $Y$.

For a square Puiseux matrix $Y\in K^{N\times N}$, let 
$$
W^{(\Lambda)}(Y)=\{\V{v}\in (R/I(\Lambda))^N\,\vert\,Y\V{v}\equiv \V{0}\mod{I(v(Y)+\Lambda)}\}
$$ 
be the set of {\it $\Lambda$-th approximate null vectors of $Y$}.
For $\lambda\in K$, we call 
$$
W^{(\Lambda)}(Y;\lambda):=W^{(\Lambda)}(Y-\lambda\cdot \mathrm{id}.)
$$ 
the set of {\it $\Lambda$-th approximate eigenvectors associated with $\lambda$}.
If $W^{(\Lambda)}(Y;\lambda)\neq\{\V{0}\}$, $\lambda$ is called a {\it $\Lambda$-th approximate eigenvalue of $Y$}.
$W^{(\Lambda)}(Y;\lambda)$ has the structure of left $R/I(\Lambda)$-module.

\begin{lemma}\label{Lemma2.1}
There exist finitely many elements $\V{v}^1,\dots,\V{v}^g\in W^{(\Lambda)}(Y;\lambda)$ such that 
\begin{equation}\label{eq1}
W^{(\Lambda)}(Y;\lambda)= \left\{(R/I(\Lambda))\cdot \V{v}^1\oplus\cdots\oplus (R/I(\Lambda))\cdot\V{v}^r\right\}
\oplus 
\left\{(I/I(\Lambda))\cdot\V{v}^{r+1}\oplus\cdots\oplus (I/I(\Lambda))\cdot\V{v}^g\right\}.
\end{equation}
These vectors satisfy $1\leq i\leq r\Rightarrow v(\V{v}^i)=0$.
\end{lemma}
\proof

For sufficiently large $d,v\in\NN$, all entries of $(Y-\lambda\cdot\mathrm{id}.)$ are contained in $t^{-v}\cdot \CC[[t^{1/d}]]$.
Because $\CC[[t^{1/d}]]$ is a PID, there exist invertible matrices $P,Q$ such that
$$
Y-\lambda\cdot\mathrm{id}.=P\cdot 
\mathrm{diag}(
\overbrace{t^{m_{1}},\dots,t^{m_r}}^{m_\ast > \Lambda+v},
\overbrace{t^{m_{r+1}},\dots,t^{m_g}}^{v<m_\ast \leq  \Lambda+v},
{t^{v},\dots,t^{v}}
)
\cdot Q,\quad (+\infty\geq m_{1}\geq \dots\geq m_{g}).
$$ 
(We regard $t^{+\infty}=0$).
It is enough to define $\V{v}^i:=Q^{-1}\V{e}_i$, where $\V{e}_i\in \RR^N$ is the $i$-th fundamental vector.
$\qed$

\begin{defi}
In the situation of (\ref{eq1}), we denote $W^{(\Lambda)}(Y;\lambda)=\kagi{\V{v}^1,\dots,\V{v}^r,t^\emptyset\V{v}^{r+1},\dots,t^\emptyset\V{v}^g}$, where $t^\emptyset$ means ``an element of infinitesimally small valuation".
\end{defi}

Although the set $\{\V{v}^1,\dots,\V{v}^g\}$ is not always a $R/I(\Lambda)$-basis of $W^{(\Lambda)}(Y;\lambda)$, it should sound less absurd if we use the term ``quasi-basis" to refer it.

\begin{defi}
We say {\it $\{\V{v}^1,\dots,\V{v}^r,t^\emptyset\V{v}^{r+1},\dots,t^\emptyset\V{v}^g\}$ to be a quasi basis of} $W^{(\Lambda)}(Y;\lambda)$ if the following conditions hold:
\begin{enumerate}
\item $W^{(\Lambda)}(Y;\lambda)=\kagi{\V{v}^1,\dots,\V{v}^r,t^\emptyset\V{v}^{r+1},\dots,t^\emptyset\V{v}^g}$,
\item $\tp{\V{v}^1},\dots,\tp{\V{v}^g}$ are linearly independent.
\end{enumerate}
\end{defi}

Let $\{\V{v}^1,\dots,\V{v}^r,t^\emptyset\V{v}^{r+1},\dots,t^\emptyset\V{v}^g\}$ be a quasi basis of $W^{(\Lambda)}(Y;\lambda)$.
From the proof of Lemma \ref{Lemma2.1}, we can soon conclude that the numbers $r$ and $g$ do not depend on the choice of quasi basis.

\begin{lemma}\label{kernel}
$\mathrm{Ker}(Y_v-\lambda_v\cdot \mathrm{id.})=\CC\tp{\V{v}^1}\oplus\cdots\oplus \CC\tp{\V{v}^g}$, where $v=v(Y)$.
\end{lemma}
\proof Let $\tilde{v}:=v(Y-\lambda\cdot\mathrm{Id})$.
By definition, we have $(Y-\lambda\cdot\mathrm{id}.)\V{v}^i\equiv \V{0}\mod{I[\tilde{v}+\Lambda]}$.
Comparing the leading terms on both sides, we obtain $(Y_v-\lambda_v\cdot\mathrm{id}.)\tp{\V{v}^i}=\V{0}$, which implies $\tp{\V{v}^i}\in \mathrm{Ker}(Y_v-\lambda_v\cdot \mathrm{id.})$.
Next assume $\V{v}\in \mathrm{Ker}(Y_v-\lambda_v\cdot \mathrm{id.})$.
Then, we have $(Y-\lambda\cdot\mathrm{id.})\V{v}\equiv \V{0}\mod{I(\tilde{v})}$.
Equivalently, $(Y-\lambda\cdot\mathrm{id.})(t^{\Lambda}\V{v})\equiv \V{0}\mod{I(\tilde{v}+\Lambda)}$.
By (\ref{eq1}), we have $t^\Lambda\V{v}\equiv r_1\V{v}^1+\cdots+r_g\V{v}^g\mod{I(\Lambda)}$ for some $r_1,\dots,r_g\in R$.
Comparing the leading terms on both sides, we obtain $\V{v}\in \CC\tp{\V{v}^1}\oplus\cdots\oplus \CC\tp{\V{v}^g}$.
\qed

\subsubsection{Restriction to $\CC((t^{1/q}))$}

Let $M=\kagi{\V{v}^1,\dots,\V{v}^r,t^\emptyset\V{v}^{r+1},\dots,t^\emptyset\V{v}^g}$ be a $R/I(\Lambda)$-module.
Set $L_q:=\CC[[t^{1/q}]]/(\CC[[t^{1/q}]]\cap I(\Lambda))$.
From the inclusion $L_q\hookrightarrow R/I(\Lambda)$, the module $M$ can be regarded as a $L_q$-module by restriction.
Especially, $M\cap L_q^N$ has a structure of  $L_q$-module.
Although $M$ is not finitely generated in general, $M\cap L_q^N$ is always finitely generated as a $L_q$-module for sufficiently large $q$.
We have
\begin{equation}\label{zure}
M\cap L_q^N=\left\{L_q\cdot \V{v}^1\oplus\cdots\oplus L_q\cdot\V{v}^r\right\}\oplus \left\{L_q\cdot (t^{1/q}\V{v}^{r+1})\oplus\cdots\oplus L_q\cdot(t^{1/q}\V{v}^g)\right\}.
\end{equation}

\subsubsection{Transposed matrix}

Let ${}^TY$ be the transposed matrix of $Y$.
Note the equation $\mathrm{rk}W^{(\Lambda)}(Y;\lambda)=\mathrm{rk}W^{(\Lambda)}({}^TY;\lambda)$, which follows from
the proof of Lemma \ref{Lemma2.1}.
Let 
$$
W^{(\Lambda)}(Y;\lambda)=\kagi{\V{x}^1,\dots,\V{x}^r,t^\emptyset\V{x}^{r+1},\dots,t^\emptyset\V{x}^g},\qquad W^{(\Lambda)}({}^TY;\lambda)=\kagi{\V{y}^1,\dots,\V{y}^r,t^\emptyset\V{y}^{r+1},\dots,t^\emptyset\V{y}^g}.
$$
As all Puiseux vectors $\V{x}^1,\dots,\V{x}^g,\V{y}^1,\dots,\V{y}^g$ are contained in $L_q^N$ for sufficiently large $q$, they expand as
\begin{equation}\label{Zure}
\begin{gathered}
\V{x}^i=\V{x}^i_{0}+\V{x}^i_{1/q}t^{1/q}+\V{x}^i_{2/q}t^{2/q}+\cdots,\quad \V{y}^i=\V{y}^i_{0}+\V{y}^i_{1/q}t^{1/q}+\V{y}^i_{2/q}t^{2/q}+\cdots,\quad (1\leq i\leq r)\\
t^{1/q}\V{x}^i=\V{x}^i_{1/q}t^{1/q}+\V{x}^i_{2/q}t^{2/q}+\cdots,\quad t^{1/q}\V{y}^i=\V{y}^i_{1/q}t^{1/q}+\V{y}^i_{2/q}t^{2/q}+\cdots,\quad (r< i\leq g).
\end{gathered}
\end{equation}
\begin{lemma}
Assume $\V{x}^i,\V{y}^i\in L^N_q$ and
$v(Y-\lambda\cdot \mathrm{id})+\Lambda=k/q$.
Then, for each $i,j$,
\begin{equation}\label{eq3a}
\sum_{s+u=k+1}{}^T\V{y}^i_0{(Y_{\frac{s}{q}}-\lambda_{\frac{s}{q}})}\V{x}^j_{\frac{u}{q}}=\sum_{s+u=k+1}{}^T\V{y}^i_{\frac{u}{q}}{(Y_{\frac{s}{q}}-\lambda_{\frac{s}{q}})}\V{x}^j_{0}.
\end{equation}
\end{lemma}
\proof In this proof, we simply write 
$\V{x}=
\left\{
\begin{array}{cc}
\V{x}^j & (1\leq j\leq r)\\
\V{x}^jt^{1/q} & (r< j\leq g)
\end{array}
\right.
$, 
$\V{y}=
\left\{
\begin{array}{cc}
\V{y}^i & (1\leq i\leq r)\\
\V{y}^it^{1/q} & (r< i\leq g)
\end{array}
\right.
$.
Because $\V{x}\in W^{(\Lambda)}(Y;\lambda)$ and $\V{y}\in W^{(\Lambda)}({}^TY;\lambda)$, we have 
\begin{align*}
{}^T\V{y}(Y-\lambda)\V{x}
&={}^T\V{y}\cdot \left(\sum_{s+u=k+1}(Y_{\frac{s}{q}}-\lambda_{\frac{s}{q}})\V{x}_{\frac{u}{q}}t^{\Lambda+v+\frac{1}{q}}+o(t^{\Lambda+v+\frac{1}{q}})\right)\\
&= \left(\sum_{s+u=k+1}\V{y}_{\frac{u}{q}}(Y_{\frac{s}{q}}-\lambda_{\frac{s}{q}})t^{\Lambda+v+\frac{1}{q}}+o(t^{\Lambda+v+\frac{1}{q}})\right)\V{x}.
\end{align*}
Comparing the coefficients of $t^{\Lambda+v+1/q}$ on both sides, we find the desired equation. 
$\qed$

\subsubsection{Eigen-quadruple}

Consider a quadruple
\begin{equation}\label{quint}
\mathcal{X}:=(Y;\lambda;\{\V{x}^1,\dots,\V{x}^r,t^\emptyset\V{x}^{r+1},\dots,t^\emptyset\V{x}^g\},\{\V{y}^1,\dots,\V{y}^r,t^\emptyset\V{y}^{r+1},\dots,t^\emptyset\V{y}^g\}),\quad
\begin{array}{l}
Y\in K^{N\times N},\\
\lambda\in K,\\
\V{x}^i,\V{y}^i\in (R/I(\Lambda))^N.
\end{array}
\end{equation}
\begin{defi}\label{EQ}
A quadruple (\ref{quint}) is called an {\it eigen-quadruple} (or an {\it EQ}, shortly) {\it of depth $\Lambda$} if it satisfies:
\begin{enumerate}
\item
$\{\V{x}^1,\dots,\V{x}^r,t^\emptyset\V{x}^{r+1},\dots,t^\emptyset\V{x}^g\}$ is a quasi basis of $W^{(\Lambda)}(Y;\lambda)$.
\item
$\{\V{y}^1,\dots,\V{y}^r,t^\emptyset\V{y}^{r+1},\dots,t^\emptyset\V{y}^g\}$ is a quasi basis of $W^{(\Lambda)}({}^TY;\lambda)$.
\end{enumerate}
\end{defi}
\begin{lemma}
Let (\ref{quint}) be an EQ. Then, $v(\V{x}^i)=v(\V{y}^i)$ for $i=1,\dots,g$.
\end{lemma}
\proof It follows from the proof of Lemma \ref{Lemma2.1}. $\qed$

For an EQ $\mathcal{X}$ (\ref{quint}), we call the number $r$ in Definition \ref{EQ} the {\it rank of $\mathcal{X}$}, and call the number $g$ 
the {\it size of $\mathcal{X}$}.
We denote $\mathrm{rk}(\mathcal{X}):=r$ and $\mathrm{sz}(\mathcal{X}):=g$. 

Define two matrices $P(\mathcal{X})=(P_{i,j})$, $Q(\mathcal{X})=(Q_{i,j})$ of size $g$ by
\[
P_{i,j}:={}^T\tp{\V{y}^i}\cdot \left(Y_{\Lambda+v+\frac{1}{q}}\V{x}_0^j +\sum_{k=1}^{\Lambda q}(Y_{\Lambda+v+\frac{1}{q}-\frac{k}{q}}-\lambda_{\Lambda+v+\frac{1}{q}-\frac{k}{q}})\V{x}^j_{\frac{k}{q}}\right),\quad
Q_{i,j}={}^T\tp{\V{y}^i}\cdot\V{x}^j_0,
\]
where $\V{x}^i_j$ and $\V{y}^i_j$ are the complex vectors in (\ref{Zure}).
Because (i) $j>\mathrm{rk}(\mathcal{X})$ implies $\V{x}^j_0=0$, (ii) $i>\mathrm{rk}(\mathcal{X})$ implies $\V{y}^i_0=0$ and (iii) $i\leq \mathrm{rk}(\mathcal{X})$ implies $\tp{\V{y}^i}=\V{y}^i_0$, the matrices $P(\mathcal{X})$ and $Q(\mathcal{X})$ are of the form
\begin{equation}\label{eq9}
P(\mathcal{X})=\left( 
\begin{array}{c|ccc}
  \Delta &  & O &  \\ \hline
   &  &  &  \\ 
  A &  & \Gamma &  \\ 
   &  &  & 
\end{array} 
\right) ,\
Q(\mathcal{X})=\left( 
\begin{array}{c|ccc}
  \Omega &  & O &  \\ \hline
   &  &  &  \\ 
  B &  & O &  \\ 
   &  &  & 
\end{array} 
\right) ,\qquad
\begin{array}{cl}
\Delta,\Omega& \in \CC^{\mathrm{rk}(\mathcal{X})\times \mathrm{rk}(\mathcal{X})},\\
A,B & \in \CC^{(\mathrm{sz}(\mathcal{X})- \mathrm{rk}(\mathcal{X}))\times \mathrm{rk}(\mathcal{X})},\\
\Gamma & \in \CC^{(\mathrm{sz}(\mathcal{X})- \mathrm{rk}(\mathcal{X}))\times (\mathrm{sz}(\mathcal{X})- \mathrm{rk}(\mathcal{X}))}
\end{array}
.
\end{equation}
We call the submatrices $\Delta=\Delta(\mathcal{X})$ and $\Omega=\Omega(\mathcal{X})$ the {\it principal matrix} of $P(\mathcal{X})$ and $Q(\mathcal{X})$, respectively.

\begin{example}
Let $Y=\left(
\begin{array}{cc}
2+t & t^2 \\ 
-t & 2+2t
\end{array} \right)$.
We have $W^{(0)}(Y;\lambda)=\kagi{\V{e}_1,\V{e}_2}$ if $\lambda=2+o(t^0)$, and otherwise $W^{(0)}(Y;\lambda)=\{\V{0}\}$.
Let $\mathcal{X}_0:=(Y;2;\{\V{e}_1,\V{e}_2\},\{\V{e}_1,\V{e}_2\})$ be a EQ of depth $0$.
It follows that $\mathrm{rk}(\mathcal{X}_0)=\mathrm{sz}(\mathcal{X}_0)=2$ and
$
P(\mathcal{X}_0)=\left(
\begin{array}{cc}
1 & 0 \\ 
-1 & 2
\end{array} 
\right)$ and 
$
Q(\mathcal{X}_0)=\left(
\begin{array}{cc}
1 & 0 \\ 
0 & 1
\end{array} 
\right)
$.
\end{example}

If $\mathcal{X}=(Y;\lambda;\{\V{x}^i\},\{\V{y}^i\})$ is an EQ of depth $\Lambda$, the quadruple ${}^T\mathcal{X}=({}^TY;\lambda;\{\V{y}^i\},\{\V{x}^i\})$ is also an EQ of depth $\Lambda$.

\begin{lemma}
$\Delta({}^T\mathcal{X})={}^T\Delta(\mathcal{X})$, $\Omega({}^T\mathcal{X})={}^T\Omega(\mathcal{X})$.
\end{lemma}
\proof Let $P({}^T\mathcal{X})=(\hat{P}_{i,j})$ and $Q({}^T\mathcal{X})=(\hat{Q}_{i,j})$.
If $1\leq i\leq \mathrm{rk}(\mathcal{X})=\mathrm{rk}({}^T\mathcal{X})$, then $\tp{\V{x}^i}=\V{x}^i_0$ and $\tp{\V{y}^i}=\V{y}^i_0$ by definition.
Therefore, from (\ref{eq3a}), we have 
\begin{align*}
P_{i,j} & ={}^T\V{y}^i_0 Y_{\Lambda+v+\frac{1}{q}}\V{x}_0^j 
+
{}^T\V{y}^i_0\cdot\sum_{k=1}^{\Lambda q}(Y_{\Lambda+v+\frac{1}{q}-\frac{k}{q}}-\lambda_{\Lambda+v+\frac{1}{q}-\frac{k}{q}})\V{x}^j_{\frac{k}{q}}\\
& ={}^T\V{y}^i_0 Y_{\Lambda+v+\frac{1}{q}}\V{x}_0^j +\sum_{k=1}^{\Lambda q}{}^T\V{y}^i_{\frac{k}{q}}(Y_{\Lambda+v+\frac{1}{q}-\frac{k}{q}}-\lambda_{\Lambda+v+\frac{1}{q}-\frac{k}{q}})\V{x}^j_{0}=\hat{P}_{j,i}
\end{align*}
and $Q_{i,j}={}^T\V{y}_0^i\V{x}^j_0=\hat{Q}_{j,i}$ for $1\leq i,j\leq \mathrm{rk}(\mathcal{X})$. $\qed$

\subsection{Simple eigenvalues}\label{sec3.1.1}

We start with the simplest case where $\tp{Y}=Y_v$ has a simple eigenvalue $\lambda_v$.
Let $q:=q(Y)$, $s:=t^{1/q}$, $Z_n:=Y_{n/q}$, $\mu_n:=\lambda_{n/q}$ and $V:=vq$.

By assumption, there exist two vectors $\V{x}_0,\V{w}\in\CC^N$ such that $(Z_V-\mu_V)\V{x}_0=\V{0}$, $({}^TZ_V-\mu_V)\V{w}=\V{0}$, ${}^T\V{w}\cdot \V{x}_0=1$
and
${}^T\V{w}\cdot \V{z}=0\iff \V{z}\in \mathrm{Im}(Z_V-\mu_V)$, ($\V{z}\in \CC^N$).

Let $f:\mathrm{Im}(Z_V-\mu_V)\to \CC^N$ be a linear map which satisfies $(Z_V-\mu_V)\circ f=\mathrm{id.}\vert_{\mathrm{Im}(Z_V-\mu_V)}$.
We recursively define countably many complex vectors $\V{x}_{1},\V{x}_{2},\dots$ and complex numbers $\mu_{V+1},\mu_{V+2},\dots$ by the following formulas (these objects are determined as: $\mu_{V+1}\to \V{x}_{1}\to \mu_{V+2}\to \V{x}_{2}\to\cdots$)
\begin{gather}\label{eq6}
\mu_{V+n}:={}^T\V{w}\cdot\left(Z_{V+n}\V{x}_0+
\sum_{i=1}^{n-1}(Z_{V+i}-\mu_{V+i})\V{x}_{n-i}\right),\quad
\V{x}_{n}:=-f\left(
\sum_{i=1}^{n}
(Z_{V+i}-\mu_{V+i})\V{x}_{n-i}
\right). 
\end{gather}

We should prove that this algorithm is well-worked.
\begin{lemma}
$\sum_{i=1}^{n}(Z_{V+i}-\mu_{V+i})\V{x}_{n-i}\in\im(Z_V-\mu_V)$.
\end{lemma}
\proof 
It suffices to prove $\sum_{i=1}^{n}{}^T\V{w}(Z_{V+i}-\mu_{V+i})\V{x}_{n-i}=\V{0}$.
We prove it by induction on $n\geq 1$.
If $n=1$, we have ${}^T\V{w}(Z_1-\mu_1)\V{x}_0=
{}^T\V{w}(Z_1-{}^T\V{w}Z_1\V{x}_0)\V{x}_0=\V{0}$.
Let $n>1$ and assume that $\V{x}_1,\dots,\V{x}_{n-1}$ are already defined. 
We have
\begin{align*}
&\sum_{i=1}^{n}{}^T\V{w}(Z_i-\mu_i)\V{x}_{n-i}
=
\sum_{i=1}^{n-1}{}^T\V{w}(Z_i-\mu_i)\V{x}_{n-i}
+{}^T\V{w}(Z_n-\mu_n)\V{x}_{0}\\
&=
\sum_{i=1}^{n-1}{}^T\V{w}(Z_i-\mu_i)\V{x}_{n-i}
+{}^T\V{w}\left(Z_n-({}^T\V{w}Z_n\V{x}_0+\sum_{i=1}^{n-1}
{}^T\V{w}(Z_i-\mu_i)\V{x}_{n-i}
)\right)\V{x}_0=\V{0}.
\end{align*}
The second equality is derived from the definition of $\mu_n$.
$\qed$

Let $\V{x}:=\sum_{i=0}^{\infty}\V{x}_{i}s^{i}$.
By direct calculations, we obtain the equation
\begin{align*}
&\sum_{n=0}^{\infty}(Z_{V+n}-\mu_{V+n})s^{V+n}\V{x}
=\sum_{n=0}^{\infty}
\sum_{i=0}^{n}(Z_{V+i}-\mu_{V+i})\V{x}_{n-i}s^{V+n}\\
&\ \ \ \ \ \ \ \ \ \ \ \ \
=\sum_{n=0}^{\infty}
\left(
(Z_V-\mu_V)\V{x}_{n}+
\sum_{i=1}^{n}(Z_{V+i}-\mu_{V+i})\V{x}_{n-i}
\right)
s^{V+n}=\V{0}.
\end{align*}
(The third equality is derived from the definition of $\V{x}_{n}$).
This implies $Z\V{x}=\mu\V{x}$, where $\mu=\sum_{i=0}^{\infty}\mu_{V+i}s^{V+i}$.

\subsection{Semi-simple eigenvalues}

Next consider the case when $\lambda_v=\mu_V$ is an semi-simple eigenvalue of $Y_v=Z_{V}$.
We note the relation
$\mathrm{Ker}({}^TZ_V-\mu_V)\cdot \V{z}=0\  \iff
\V{z}\in\mathrm{Im}(Z_V-\mu_V)$, ($\V{z}\in \CC^N$).


\begin{lemma}\label{lemma3.8a}
Let $\mathcal{X}=(Z;\mu;\{\V{x}^i\},\{\V{y}^i\})$ be an EQ of depth $\Lambda$, rank $r$ and size $g$.
The submatrix $\Gamma$ of $P(\mathcal{X})$ in (\ref{eq9}) is invertible. 
\end{lemma}
\proof 
If $\Gamma$ is not invertible, there must exist some nonzero vector $(0,\dots,0,c_{p},c_{p+1},\dots,c_g)\in \CC^g$ ($r<p\leq g$, $c_p\neq 0$) such that the vector $\V{z}:=\sum_{i=p}^{g}c_i\V{x}^i=\V{z}_1s+\V{z}_2s^2+\cdots+\V{z}_\Lambda s^\Lambda$ satisfies 
$$
\tp{{}^T\V{y}^i}\cdot \sum_{k=1}^\Lambda(Z_{\Lambda+V+1-k}-\mu_{\Lambda+V+1-k})\V{z}_{k}=\V{0}
$$ 
for $i=1,\dots,g$.
From Lemma \ref{kernel}, this is equivalent to
$
\mathrm{Ker}({}^TZ_V-\mu_V)\cdot \sum (Z_{\Lambda+V+1-k}-\mu_{\Lambda+V+1-k})\V{z}_{k}=\V{0}
$.
Therefore, there exists some complex vector $\V{v}$ such that 
\begin{equation}\label{eq8}
(Z_V-\mu_V)\V{v}+\sum_{k=1}^\Lambda(Z_{\Lambda+V+1-k}-\mu_{\Lambda+V+1-k})\V{z}_{k}=\V{0}.
\end{equation}
Set $\V{z}^\ast:=\V{z}_{1}+\V{z}_{2}s+\cdots+\V{z}_{\Lambda}s^{\Lambda-1}+\V{v}s^{\Lambda}$.
By $(Z-\mu)\V{z}\equiv \V{0}\mod{I(\Lambda+V)}$ and (\ref{eq8}), it is directly checked that $(Z-\mu)\V{z}^\ast\equiv 0\mod{I(\Lambda+V)}$. 
Therefore, $\V{z}^\ast\in W^{(\Lambda)}(Z;\mu)$.
By construction, we have  $\tp{\V{z}^\ast}=c_{p}\tp{\V{x}^{p}}+\dots+c_{g}\tp{\V{x}^{g}}$.
On the other hand, because $\{\V{x}^i\}$ is an quasi basis of $W^{(\Lambda)}(Z;\mu)$, we have $\V{z}^\ast=a_1\V{x}^1+\dots+a_{g}\V{x}^{g}$ for some $a_1,\dots,a_{g}\in R$.
As $1\leq i\leq r\iff v(\V{x}^i)=0$, the equation $\tp{\V{z}^\ast}=\tp{(a_1\V{x}^1)}+\dots+\tp{(a_r\V{x}^r)}$ holds.
Comparing these two equations, we obtain $c_{p}\tp{\V{x}^{p}}+\dots+c_{s}\tp{\V{x}^{s}}=\tp{(a_1\V{x}^1)}+\dots+\tp{(a_{r}\V{x}^{r})}$, which contradicts the fact that $\{\tp{\V{x}^i}\}$ is linearly independent. $\qed$

\begin{cor}\label{cor2.13}
Let $\zeta\in \CC$ be a complex number.
The correspondence
\[
\left\{
\V{c}=\left( 
\begin{array}{c}
\V{v} \\ 
\V{v}'
\end{array} 
\right) \in\CC^g\,;\,\{P(\mathcal{X})-\zeta\cdot Q(\mathcal{X})\}\V{c}=\V{0}
\right\}
\to 
\left\{
\V{v}\in\CC^r\,;\,\{\Delta(\mathcal{X})-\zeta\cdot\Omega(\mathcal{X})\}\V{v}=\V{0}
\right\},
\quad
\left( 
\begin{array}{c}
\V{v} \\ 
\V{v}'
\end{array} 
\right) 
\mapsto \V{v}
\]
is a linear isomorphism.
\end{cor}

\proof It is straightforward from Lemma \ref{lemma3.8a} and (\ref{eq9}).
$\qed$

Let 
$L(\zeta):=\{\V{c}\in\CC^g\,\vert\,\{P(\mathcal{X})-\zeta\cdot Q(\mathcal{X})\}\V{c}=\V{0} \}$ and $\{\V{c}^1,\dots,\V{c}^{\rho}\}$ be a basis of $L(\zeta)$ ($0\leq \rho=\dim_\CC L(\zeta)\leq r$).
Put $(\V{X}^1,\dots,\V{X}^\rho):=(\V{x}^1,\dots,\V{x}^g)\cdot (\V{c}^1,\dots,\V{c}^\rho)$.
The vector $\V{X}^i=\V{X}^i_0+\V{X}_{1}^is^{1}+\cdots+\V{X}_\Lambda s^\Lambda$ satisfies the equation
$$
{}^T\tp{\V{y}^i}(Z_{\Lambda+V+1}-\zeta)\V{X}^j_0+\sum_{k=1}^\Lambda{}^T\tp{\V{y}^i}(Z_{\Lambda+V+1-k}-\mu_{\Lambda+V+1-k})\V{X}^j_{k}=\V{0}
$$ 
($\forall i,j$), which is equivalent to 
$$
\mathrm{Ker}{({}^TZ_V-\mu_V)}\cdot
\left\{
(Z_{\Lambda+V+1}-\zeta)\V{X}^j_0+\sum
(Z_{\Lambda+V+1-k}-\mu_{\Lambda+V+1-k})\V{X}^j_{k}
\right\}=\V{0}.
$$
Therefore, there exists some complex vector $\V{X}^j_{\Lambda+1}$\footnote{The vector $\V{X}^j_{\Lambda+1}$ is calculated by the formula 
\[
\V{X}^j_{\Lambda+1}=-f((Z_{\Lambda+V+1}-\zeta)\V{X}^j_0+(Z_{\Lambda+V}-\mu_{\Lambda+V})\V{X}^j_{1}+\cdots+(Z_{V+1}-\mu_{V+1})\V{X}^j_{\Lambda}),
\]
where $f:\mathrm{Im}(Z_V-\mu_V)\to \CC^N$ is a linear map satisfying $(Z_V-\mu_V)\circ f=\mathrm{Id}$. 
} such that 
\[
(Z_{\Lambda+V+1}-\zeta)\V{X}^j_0+\sum_{k=1}^\Lambda(Z_{\Lambda+V+1-k}-\mu_{\Lambda+V+1-k})\V{X}^j_{k}
+(Z_V-\mu_V)\V{X}^j_{\Lambda+1}=\V{0}.
\]
This means that the vector $\tilde{\V{x}}^j:=\V{X}^j_0+\V{X}^j_{1}s^{1}+\dots+\V{X}^j_{\Lambda+1}s^{\Lambda+1}$ and an element $\tilde{\mu}:=\mu_Vs^V+\mu_{V+1}s^{V+1}+\dots+\mu_{V+\Lambda}s^{V+\Lambda}+\zeta s^{V+\Lambda+1}$ satisfy $\tilde{\V{x}}^j\in W^{(\Lambda+1)}(Z;\tilde{\mu})$.

\begin{algo}\label{newEQ}
We define the new EQ
\begin{equation}
\tilde{\mathcal{X}}=(Z;\tilde{\mu};\{\tilde{\V{x}}^i\},\{\tilde{\V{y}}^i\})
\end{equation}
of depth $\Lambda+1$ by the following manner:
\begin{enumerate}
\renewcommand{\theenumi}{(\arabic{enumi})}
\item Define $\tilde{\mu},\tilde{\V{x}}^1,\dots,\tilde{\V{x}}^\rho$ as above.
\item Pick up $\tilde{\V{x}}^{\rho+1},\dots,\tilde{\V{x}}^{r}\in \{\V{x}^1,\dots,\V{x}^r\}$ so that $\tp{\tilde{\V{x}}^1},\dots,\tp{\tilde{\V{x}}^r}$ are linearly independent.
\item Define $\tilde{\V{x}}^{i}:=\V{x}^{i}\cdot s$ for $r<i\leq g$.
\item Then, we obtain the new quasi basis  $W^{(\Lambda+1)}(Y;\tilde{\lambda})=\kagi{\tilde{\V{x}}^1,\dots,\tilde{\V{x}}^\rho,s^\emptyset\tilde{\V{x}}^{\rho+1},\dots,s^\emptyset\tilde{\V{x}}^g}$.
\item The Puiseux vectors $\tilde{\V{y}}^i$ can be defined by similar manner.\qed
\end{enumerate}
\end{algo}

\begin{rem}
We must note the fact that the matrix equation $\{P(\mathcal{X})-\zeta\cdot Q(\mathcal{X})\}\V{c}=\V{0}$ does not always have a solution $(\zeta,\V{c})$, ($\V{c}\neq \V{0}$).
This is a typical case of failure where the naive algorithm is not applicable. 
(See Examples \ref{ex2.17} and Remark \ref{ex2.18}).
We will avoid this difficulty by conjugating the matrix $Y$ by some diagonal matrix. 
\end{rem}

\subsection{The case if $\mathrm{rk}(\mathcal{X})=1$}\label{rk1}

If $\mathrm{rk}(\mathcal{X})=1$, we can make our algorithm simpler.
Let $\mathcal{X}=(Z;\mu;\{\V{x}^1,s^\emptyset\V{x}^2,\dots,s^\emptyset\V{x}^g \},\{\V{y}^1,s^\emptyset\V{y}^2,\dots,s^\emptyset\V{y}^g \})$ be an EQ of rank $1$.
Let $f:\mathrm{Im}(Z_v-\mu_V)\to \CC^N$ be a linear map such that $(Z_V-\mu_V)\circ f=\mathrm{Id}$.

Because $\mu_V$ is semi-simple, we can assume ${}^T\V{y}_0^1\cdot \V{x}^1_0=1$ without loss of generality.
The matrices $P(\mathcal{X})$, $Q(\mathcal{X})$ are of the form:
\[
P(\mathcal{X})=\left(
\begin{array}{c|ccc}
P_{1,1} &  & O &  \\ \hline
 &  &  &  \\ 
\V{p} &  & \Gamma &  \\ 
 &  &  & 
\end{array} 
\right),\qquad
Q(\mathcal{X})=\left(
\begin{array}{c|ccc}
1 &  & O &  \\ \hline
 &  &  &  \\ 
\V{q} &  & O &  \\ 
 &  &  & 
\end{array} 
\right)
\]
Obviously, the matrix equation $\{P(\mathcal{X})-\zeta\cdot Q(\mathcal{X})\}\V{c}=\V{0}$ has a solution $\zeta=P_{1,1}$ and $\V{c}={}^T(1,\V{\delta})$, where $\V{\delta}=\Gamma^{-1}(P_{1,1}\V{q}-\V{p})$.
Set
\begin{gather*}
\tilde{\V{x}}^1 %
:=\V{x}^1+(\V{x}^2s,\dots,\V{x}^gs)\cdot \V{\delta}-f\left((Z_{\Lambda+V+1}-\zeta)\tilde{\V{x}}^1_0+(Z_{\Lambda+V}-\mu_{\Lambda+V})\tilde{\V{x}}^1_{1}+\cdots+(Z_{V+1}-\mu_{V+1})\tilde{\V{x}}^1_\Lambda\right),
\end{gather*}
where $\V{x}^1+(\V{x}^2s,\dots,\V{x}^gs)\cdot \V{\delta}=\tilde{\V{x}}^1_0+\tilde{\V{x}}^1_{1}s^{1}+\cdots+\tilde{\V{x}}^1_{\Lambda}s^{\Lambda}$.
Therefore, the quintuple 
\[
\tilde{\mathcal{X}}=(Z;\mu+\zeta s^{\Lambda+1};\{\tilde{\V{x}}^1,s^{1+\emptyset}\V{x}^2,\dots,s^{1+\emptyset}\V{x}^g \},\{\tilde{\V{y}}^1,s^{1+\emptyset}\V{y}^2,\dots,s^{1+\emptyset}\V{y}^g \})
\]
is an EQ whose depth is greater then that of $\mathcal{X}$. 
Obviously, we have ${}^T\tp{\tilde{\V{x}}^1}\cdot \tp{\tilde{\V{y}}^1}={}^T\V{x}^1_0\cdot \V{y}^1_0=1$, which admits us to repeat this procedure.

\subsection{Examples}

\begin{example}[Simple eigenvalue]
Let $Y=
\left(
\begin{array}{ccc}
1+t^2 & 2t  & 2 \\ 
1+t & 2-t & 2t \\ 
2 & t^2 & 1+t^2
\end{array} 
\right)=:
Y_0+Y_1t+Y_2t^2
$.
The leading term of $Y$ is 
$Y_0=\left(
\begin{array}{ccc}
1 & 0  & 2 \\ 
1 & 2 & 0 \\ 
2 & 0 & 1
\end{array} 
\right)$, which has a simple eigenvalue $\lambda_0=4$ and its corresponding eigenvector $\V{x}_0
=
\left(\begin{array}{c}
1 \\ 
1 \\ 
1
\end{array} \right)
$.
The transposed matrix ${}^TY_0$ has also an eigenvector 
$\V{w}=
\dfrac{1}{2}\left(\begin{array}{c}
1 \\ 
0 \\ 
1
\end{array} \right)$ corresponding to $4$.
Fix a linear map $f:\mathrm{Im}(Y_0-4\mathrm{Id}.)\to \CC^3$ which satisfies $(Y_0-4\mathrm{Id}.)\circ f=\mathrm{Id.}_{\mathrm{Im}(Y_0-4\mathrm{Id}.)}$. For example, define $f(\V{e}_1-\V{e}_3):=-\frac{1}{2}\V{e}_1-\frac{1}{2}\V{e}_2$ and $f(\V{e}_2):=-\V{e}_2$.

From (\ref{eq6}), we have $\lambda_1={}^T\V{w}Y_1\V{x}_0=1$, $\V{x}_1=-f((Y_1-1)\V{x}_0)=\tfrac{1}{2}\V{e}_1+\tfrac{3}{2}\V{e}_2$, $\lambda_2={}^T\V{w}\{Y_2\V{x}_0+(Y_0-1)\V{x}_1\}=\tfrac{9}{4}$, $\V{x}_2=-f((Y_2-\tfrac{9}{4})\V{x}_0+(Y_1-1)\V{x}_1)=\tfrac{5}{8}\V{e}_1-\tfrac{33}{8}\V{e}_2,\dots$.
Finally, we obtain an eigenvector 
$$
\V{x}=
\left(
\begin{array}{c}
1 \\ 
1 \\ 
1
\end{array} 
\right)+
\left(
\begin{array}{c}
1/2 \\ 
3/2 \\ 
0
\end{array} 
\right)t+
\left(
\begin{array}{c}
5/8 \\ 
-33/8 \\ 
0
\end{array} 
\right)t^2+\cdots
$$
corresponding to an eigenvalue $\lambda=4+t+\frac{9}{4}t^2+\cdots$.
\end{example}

\begin{example}[Semi-simple eigenvalue]\label{ex2.16}
Let 
$Y=\left(\begin{array}{ccc}
1 & 1-t & t \\ 
t & 2 & 2t \\ 
2t & t^2 & 2+t
\end{array} 
\right)=Y_0+Y_1t+Y_2t^2
$.
The leading term $Y_0$ has a semi-simple eigenvalue $\lambda_0=2$, whose corresponding eigenspace is generated by the two vectors 
$\V{x}^1={}^T(1,1,0)$ and 
$\V{x}^2={}^T(0,0,1)$.
The transposed matrix ${}^TY$ has eigenvectors $\V{y}^1={}^T(0,1,0)$ and $\V{y}^2={}^T(0,0,1)$.
Fix a linear map $f:\mathrm{Im}(Y_0-2\mathrm{Id})\to \CC^3$ such that $(Y_0-2\mathrm{Id})\circ f=\mathrm{Id}$, for example, define $f(\V{e}_1):=\V{e}_2$.
We obtain a sequence of EQs $\mathcal{X}_0,\mathcal{X}_1,\dots$ by the following step-by-step method:
\begin{itemize}
\item The quintuple $\mathcal{X}_0=(Y;2;\{\V{x}^1,\V{x}^2\},\{\V{y}^1,\V{y}^2\})$ is an EQ of depth $0$ with 
$
P(\mathcal{X}_0)=
\left(
\begin{array}{cc}
1 & 2 \\ 
2 & 1
\end{array} 
\right)$,
$Q(\mathcal{X}_0)=\left(
\begin{array}{cc}
1 & 0 \\ 
0 & 1
\end{array} 
\right)$.
The matrix equation $\{P(\mathcal{X}_0)-\zeta Q(\mathcal{X}_0)\}\V{c}=\V{0}$ has a solution $\zeta=3$ and $\V{c}={}^T(1,1)$, and the transposed equation $\{P({}^T\mathcal{X}_0)-\zeta Q({}^T\mathcal{X}_0)\}\V{d}=0$ has a solution $\zeta=3$, $\V{d}=\frac{1}{2}{}^T(1,1)$.
\item Define $\tilde{\V{x}}_0^1:=(\V{x}^1,\V{x}^2)\V{c}={}^T(1,1,1)$, $\tilde{\V{x}}_1^1:=-f((Y_1-3)\tilde{\V{x}}^1_0)={}^T(0,3,0)$ and $\tilde{\V{x}}^1:=\tilde{\V{x}}^1_0+\tilde{\V{x}}^1_1t$.
Also define $\tilde{\V{y}}^1_0:=(\V{y}^1,\V{y}^2)\V{d}=\frac{1}{2}{}^T(0,1,1)$.
Then, there exists some $\tilde{\V{y}}^1_1$ such that the quintuple 
$
\mathcal{X}_1=(Y;2+3t,\{\tilde{\V{x}}^1,t^\emptyset{\V{x}}^2\},\{\tilde{\V{y}}^1,t^\emptyset{\V{y}}^2\})$,
($\tilde{\V{y}}^1=\tilde{\V{y}}^1_0+\tilde{\V{y}}^1_1t$) 
is an EQ\footnote{We do not need to calculate $\tilde{\V{y}}^1_1$ here.} of depth $1$ with $P(\mathcal{X}_1)=
\left(
\begin{array}{cc}
-4 & 0 \\ 
1 & -2
\end{array} 
\right)$ and
$Q(\mathcal{X}_1)=
\left(
\begin{array}{cc}
1 & 0 \\ 
1 & 0
\end{array} 
\right)
$.
The matrix equation $\{P(\mathcal{X}_1)-\zeta Q(\mathcal{X}_1)\}\V{c}=0$ has a solution $\zeta=-4$ and $\V{c}={}^T(2,5)$.
\item 
Set $\hat{\V{x}}^1_0+\hat{\V{x}}^1_1t:=(\tilde{\V{x}}^1,{\V{x}}^2t)\V{c}={}^T(2,2,2)+{}^T(0,6,5)t$, $\hat{\V{x}}^1_2:=-f((Y_2+4)\hat{\V{x}}^1_0+(Y_1-3)\hat{\V{x}}^1_1)={}^T(0,-7,0)$ and $\hat{\V{x}}^1:=\hat{\V{x}}^1_0+\hat{\V{x}}^1_1t+\hat{\V{x}}^1_2t^2$.
On the other hand, there exists some $\hat{\V{y}}^1=\hat{\V{y}}^1_0+\cdots$ ($\hat{\V{y}}^1_0=\frac{1}{4}{}^T(0,1,1)$) such that the quintuple 
$
\mathcal{X}_2=(Y;2+3t-4t^2,\{\hat{\V{x}}^1,t^{1+\emptyset}{\V{x}}^2\},\{\hat{\V{y}}^1,t^{1+\emptyset}{\V{y}}^2\})
$ 
is an EQ of depth $2$.
\end{itemize}

Repeating this procedure, we obtain a sequence of EQs $\mathcal{X}_3,\mathcal{X}_4,\dots$.
In other words, we can calculate $\Lambda$-th approximate eigenvector/value for arbitrarily large $\Lambda$. 
In fact, we soon obtain an eigenvector 
\[
\V{x}=
\left(
\begin{array}{c}
2 \\ 
2 \\ 
2
\end{array} 
\right)+
\left(
\begin{array}{c}
0 \\ 
6 \\ 
5
\end{array} 
\right)t-
\left(
\begin{array}{c}
0 \\ 
7 \\ 
0
\end{array} 
\right)t^2+\cdots
\]
and its corresponding eigenvalue $\lambda=2+3t-4t^4+\cdots$.
\end{example}

\begin{example}[Case of failure]\label{ex2.17}
Let 
$Y=\left(\begin{array}{ccc}
1 & 1-t & t \\ 
t & 2 & 2t^2 \\ 
2t & t^2 & 2+t
\end{array} 
\right)=Y_0+Y_1t+Y_2t^2
$.
By similar procedures as last examples, we soon obtain an EQ $\mathcal{X}_1=(Y;2+t,\{\V{x}^1,t^\emptyset{\V{x}}^2\},\{{\V{y}}^1,t^\emptyset{\V{y}}^2\})$ of depth $1$, where ${\V{x}}^1={}^T(0,0,1)-{}^T(0,1,0)t$, ${\V{x}}^2={}^T(0,0,1)$,  ${\V{y}}^1={}^T(0,1,0)+\cdots$ and ${\V{y}}^2={}^T(0,1,0)$.
We have
$
P(\mathcal{X}_1)=
\left(
\begin{array}{cc}
1 & 0 \\ 
1 & 2
\end{array} 
\right)$ and 
$Q(\mathcal{X}_1)=
\left(
\begin{array}{cc}
0 & 0 \\ 
1 & 0
\end{array} 
\right)$.
As the matrix equation $\{P(\mathcal{X}_1)-\zeta Q(\mathcal{X}_1)\}\V{c}=0$ has no nonzero solution $\V{c}$ for any $\zeta$, we fail to obtain an EQ of depth $2$.
This is a typical case of failure.
See the following Remark \ref{ex2.18}.
\end{example}

\begin{rem}\label{ex2.18}
Let $\delta:=\mathrm{diag}(1,t^{1/2},1)$ and $Z:=\delta^{-1}Y\delta=
\left(
\begin{array}{ccc}
1 & t^{1/2}-t^{3/2} & t \\ 
t^{1/2} & 2 & 2t^{3/2} \\ 
2t & t^{5/2} & 2+t
\end{array} 
\right)=Z_0+Z_{1/2}t^{1/2}+Z_1t+\cdots+Z_{5/2}t^{5/2}
$, where $Y$ is the Puiseux matrix defined in Example \ref{ex2.17}.
Unlike the matrix $Y$, the matrix $Z$ admits the recursive procedure:
%
%
%
\begin{itemize}
\item Let $f:\mathrm{Im}(Z_0-2\mathrm{Id}.)\to \CC^3$ be the linear map defined by $f(\V{e}_1)=-\V{e}_1$.
\item  $\mathcal{X}_{1/2}=(Z;2;\{\V{x}^1,\V{x}^2\},\{\V{y}^1,\V{y}^2\})$, $\V{x}^1={}^T(t^{1/2},1,0)$, $\V{x}^2={}^T(0,0,1)$, $\V{y}^1={}^T(0,1,0)+\cdots$, $\V{y}^2={}^T(0,0,1)+\cdots$ and $P(\mathcal{X}_{1/2})=Q(\mathcal{X}_{1/2})=\mathrm{Id}$.
The matrix equation $\{P(\mathcal{X}_{1/2})-\mu Q(\mathcal{X}_{1/2})\}\V{c}=0$ has a solution $\mu=1$ and $\V{c}=$ arbitrary, and the transposed equation $\{P({}^T\mathcal{X}_{1/2})-\mu Q({}^T\mathcal{X}_{1/2})\}\V{d}=0$ also has the same solution.
\item $\mathcal{X}_1=(Z;2+t;\{\V{x}^1,\V{x}^2\},\{\V{y}^1,\V{y}^2\})$, $\V{x}^1={}^T(t^{1/2},1,0)$, $\V{x}^2={}^T(t,0,1)$, $\V{y}^1={}^T(0,1,0)+\cdots$, $\V{y}^2={}^T(0,0,1)+\cdots$ and $P(\mathcal{X}_1)=\left(
\begin{array}{cc}
0 & 3 \\ 
2 & 0
\end{array} 
\right)$, $Q(\mathcal{X}_{1})=\mathrm{Id}$.
The matrix equation $\{P(\mathcal{X}_1)-\mu Q(\mathcal{X}_1)\}\V{c}=0$ has a solution $\mu=\sqrt{6}$ and $\V{c}={}^T(\sqrt{3},\sqrt{2})$, and the transposed equation $\{P({}^T\mathcal{X}_1)-\mu Q({}^T\mathcal{X}_1)\}\V{d}=0$ has a solution $\mu=\sqrt{6}$ and $\V{d}=\frac{1}{2\sqrt{6}}{}^T(\sqrt{2},\sqrt{3})$ 
\item $\mathcal{X}_{3/2}=(Z;2+t+\sqrt{6}t^{3/2};\{{\V{x}}^1,t^\emptyset{\V{x}}^2\},
\{{\V{y}}^1,t^\emptyset\V{y}^2 \})$, $\V{x}^1={}^T(\sqrt{3}t^{1/2}+\sqrt{2}t-2\sqrt{3}t^{3/2},\sqrt{3},\sqrt{2})$, $\V{x}^2={}^T(t,0,1)$, $\V{y}^1=\frac{1}{2\sqrt{6}}{}^T(0,\sqrt{2},\sqrt{3})+\cdots$, $\V{y}^2={}^T(0,0,1)+\cdots$.
\end{itemize}

Finally, we obtain the eigenvector of $Z$:
\[
\V{x}=
\left(
\begin{array}{c}
0 \\ 
\sqrt{3} \\ 
\sqrt{2}
\end{array} 
\right)+
\left(
\begin{array}{c}
\sqrt{3} \\ 
0 \\ 
0
\end{array} 
\right)t^{1/2}+
\left(
\begin{array}{c}
\sqrt{2} \\ 
0 \\ 
0
\end{array} 
\right)t^{1}+
\left(
\begin{array}{c}
-2\sqrt{3} \\ 
0 \\ 
0
\end{array} 
\right)t^{3/2}+\cdots
\]
corresponding to $\lambda=2+t+\sqrt{6}t^{3/2}+\cdots$.
\end{rem}

\section{Short reviews on graph theory and tropical linear algebra}\label{sec3}

Combinatorial aspects of positive matrices, complex matrices, Puiseux matrices,\,{\it etc.}\,\,reflect on their {\it (weighted) adjacency graphs}.
In this section, we discuss about {\it weighted directed graphs} and their properties.

\subsection{Weighted digraphs}

For a digraph (= directed graph) $G$, we denote its node set by $\mathrm{node}(G)$, and its arc set by $\mathrm{arc}(G)$.
When there is no chance of confusion, we simply write "$i\in A$" instead of "$i\in \mathrm{node}(A)$", and "$(i\to j)\in A$" instead of "$(i\to j)\in \mathrm{arc}(A)$".

If a {\it weighted digraph} $G$ contains an arc of weight $X$ from a node $i$ to a node $j$, we simply say "$G$ contains an arc $i\xrightarrow{X} j$".
The weight of an arc $(i\to j)\in G$ is written as $w_G(i\to j)$.

A {\it (directed) path} is a sequence of finitely many arcs: $i\to i_1\to i_2\to \dots\to j$.
The {\it length} of a path is a number of arcs contained in the path, and the {\it weight} of a path is the sum of weights of the arcs.
A {\it closed path}, or a {\it loop}, is a directed path from a node to itself.
A {\it simple loop} is a loop of length $1$.

Two nodes $i,j$ are said to be {\it strongly connected} if there exist paths both from $i$ to $j$ and from $j$ to $i$. 
A {\it strongly connected digraph} is a digraph in which any two nodes are strongly connected.
For two nodes $i,j$, we denote $i\leftrightsquigarrow j$ if a digraph contains $(i\to j)$ or $(j\to i)$.
Two nodes $i,j$ are said to be {\it connected} if there exists a sequence $i\leftrightsquigarrow i_1\leftrightsquigarrow i_2 \leftrightsquigarrow\dots \leftrightsquigarrow j$.

The {\it smallest weight} of a weighted digraph $G$ is the smallest number $s(G)$ among the weights of arcs in $G$.
If $G$ contains no arc,  $s(G):=+\infty$.
A {\it truncated graph of $G$} is the subgraph $G_{\leq\mu}\subset G$ ($\mu\in\QQ$) which is gotten by removing arcs whose weight is greater than $\mu$.
The digraph $G_{\leq s(G)}$ is called the {\it leading term of $G$}.
Denote $\tp{G}:=G_{\leq s(G)}$.

\begin{defi}
A {\it homomorphism between weighted digraphs} $G$ and $G'$ is a map $f:\mathrm{node}(G)\to \mathrm{node}(G')$ which satisfies the following property: 
\[
((i\xrightarrow{W} j) \in G) \wedge (f(i)\neq f(j)) \Longrightarrow ((f(i)\xrightarrow{W'}f(j))\in G')\wedge (W\geq W')
\]
\end{defi}
\begin{rem}
If $f:G\to G'$ be a homomorphism between weighted digraphs, $f$ induces a restricted homomorphism $f_{\leq \mu}:G_{\leq \mu}\to G'_{\leq \mu}$ for any $\mu$.
\end{rem}

Let $G$ be a weighted digraph, and let $\sim$ be an equivalence relation on $\mathrm{node}(G)$.
The {\it quotient graph} $G/\sim$ is the weighted digraph defined by
$
\mathrm{node}(G/\sim):=\mathrm{node}(G)/\sim
$
and
\begin{align*}
\mathrm{arc}(G/\sim)&:=\left\{(x\xrightarrow{W} y)\,\right\vert\left.\,x\neq y,\, 
W=\min_{i\in x,\ j\in y}\{w\,\vert\,(i\xrightarrow{w}j)\in G \}<+\infty\}
\right\}.
\end{align*}
There uniquely exists a natural homomorphism $\pi:G\to G/\sim$ of weighted digraphs such that $\pi(i)=x \iff i\in x$.

Let $B\subset G$ be a subgraph.
By identifying all nodes in $B$ to a single equivalence class and regarding each node outside of $B$ as an equivalence class with one element, we define the equivalence relation $\sim_{B}$ on $\mathrm{node}(G)$.
We will simply denote $G/B:=G/\sim_B$.

\begin{example}
Let 
$
G=
\left(
\vcenter{
\xymatrix{
1 \ar@(ld,lu)^1\ar@<0.5ex>[r]^2\ar@<0.5ex>[d]^3& 2 \ar@(rd,ru)_4\ar@<0.5ex>[l]^7 \ar[d]^6\ar[dl]^5\\
3 \ar@(ld,lu)^8\ar@<0.5ex>[u]^0 & 4
}
}
\right)
$.
Define the equivalence relation $\approx$ on $\mathrm{node}(G)$ by $i\approx j\iff i-j$ is even.
Then 
$(G/\approx)
=\left(
\vcenter{
\xymatrix{
\{1,3\} \ar@<0.5ex>[r]^2& \{2,4\} \ar@<0.5ex>[l]^5
}
}
\right)$.
\end{example}

Let $S$ be a subset of $\mathrm{node}(G)$.
A {\it subgraph induced by} $S$ is a subgraph $H\subset G$ such that $\mathrm{node}(H)=S$ and $\mathrm{arc}(H)=\{(i\to j)\in \mathrm{arc}(G)\,\vert\,i,j\in S\}$.
The induced subgraph is written as $H=G[S]$.

A digraph $T$ is said to be a {\it directed tree} (or simply, {\it tree}) if $T$ contains a node $i$ such that, for each node $j\in T$, $T$ contains exactly one path from $i$ to $j$.
A {\it forest} is a disjoint union of trees.

A digraph is said to be {\it acyclic} if it contains no closed path.
Any acyclic digraph admits a semi-order $\succeq$ defined by $i\succeq j\iff$ there exists a path from $i$ to $j$. 
We say $G$ to be an {\it acyclic digraph on $x$} if $G$ contains an unique maximal node $x$ with respected to the semi-order $\succeq$.

Let $G$ be a digraph, and $E=\{x_1,\dots,x_k\}$ be a proper subset of $\mathrm{node}(G)$. 
$G$ is said to be an {\it $E$-forest} if each connected component of $G$ (i) is acyclic, (ii) is on $x_i$ for some $i$ and (iii) $E\cap G=\{x_i\}$.

\begin{example}[$E$-forest]
The following digraph is an example of $E$-forest.
The sign $\circ$ represents a node in $E$.
The graph consists of four connected components. 
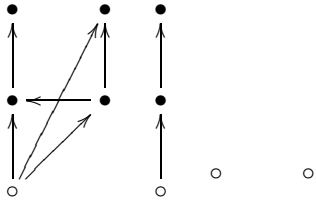
\[
\xymatrix{
\bullet & \bullet  \\
\bullet\ar[u] &\bullet\ar[l]\ar[u]  \\
\circ \ar[u]\ar[ur]\ar[uur]
}\quad
\xymatrix{
\bullet\\
\bullet\ar[u]\\
\circ \ar[u]
}\quad
\xymatrix{
 & \\
 & \\
\circ & \circ 
}
\]
\end{example}

We define the {\it eigenvalue of a weighted digraph $G$} by the formula
\begin{equation}
\Lambda(G):=\min\left\{\left.\frac{\mathrm{weight}(\gamma)}{\mathrm{length}(\gamma)}\,\right\vert\, \gamma \mbox{ is a closed path on $G$}\right\}.
\end{equation}
Of course, $s(G)\leq \Lambda(G)$.
The equality holds if and only if $\tp{G}$ contains a closed path.

\subsection{Tropical algebra}

We shortly introduce some basic facts on tropical algebra.
For details, see \cite{maclagan2015introduction}, for example.

A {\it tropical matrix} is a matrix over the {\it tropical semi-field} $\TT:=\QQ\cup \{+\infty\}$.
For two tropical matrices $C=(C_{i,j})\in \TT^{N\times M}$ and $D:=(D_{i,j})\in \TT^{M\times L}$, we define the {\it tropical product} $C\odot D\in\TT^{N\times L}$ by $C\odot D=(X_{i,j})$, $X_{i,j}:=\min_{k=1}^M[ C_{i,k}+D_{k,j}]$. 
For $a\in\TT$ and $C=(C_{i,j})\in\TT^{N\times M}$, $a+C:=(a+C_{i,j})$.

Let $C=(X_{i,j})\in \TT^{N\times N}$ be a square tropical matrix.
The {\it adjacency graph of $C$} is the weighted digraph $G(C)=(\mathrm{node}(G(C)),\mathrm{arc}(G(C)))$, where $\mathrm{node}(G(C))=\{1,\dots,N\}$ and $\mathrm{arc}(G(C))=\{j\xrightarrow{X_{i,j}} i\,\vert\,X_{i,j}\neq +\infty \}$.
A tropical matrix is said to be {\it irreducible} if its adjacency graph is strongly connected.

\begin{rem}
The $1\times 1$ matrix $(+\infty)\in \TT^{1\times 1}$ is an irreducible matrix, whose adjacency graph of is the simple node graph. 
\end{rem}

Denote $\V{+\infty}:={}^T(+\infty,\cdots,+\infty)\in \TT^N$.

\begin{prop}[Max-plus spectral theorem]
\label{trop}
Let $C\in \TT^{N\times N}$ be a tropical matrix.
\begin{enumerate}
\item There exist a tropical number $\Lambda\in \TT$ and a tropical vector $\V{V}\in\TT^N$ ($\V{V}\neq \V{+\infty}$) such that $C\odot \V{V}=\Lambda+\V{V}$.
The number $\Lambda$ is called a {\it tropical eigenvalue of $C$}, and the vector $\V{V}$ is called a {\it tropical eigenvector corresponding with $\Lambda$}.
\item If $C\neq (+\infty)$ is irreducible, its tropical eigenvalue is unique.
In this case, the tropical eigenvalue of $C$ coincides with the eigenvalue of the graph $G(C)$.
\item If $C\neq (+\infty)$, there exists a tropical eigenvector corresponding with $\Lambda$ which is contained in $\QQ^N$.
\end{enumerate}
\end{prop}
\proof See, for example, \cite{cochet1998numerical,reischuk1997} and the references therein. $\qed$

Let $\mathcal{D}(N):=\{(X_{i,j})\in \TT^{N\times N}\,\vert\,i\neq j\Leftrightarrow X_{i,j}=+\infty\}$ be the set of {\it tropical diagonal matrices of size $N$}.
We write the tropical diagonal matrix $C\in \mathcal{D}(X)$ as $C=\mathrm{diag}[T_1,\dots,T_N]$, where $T_i$ is the $(i,i)$-entry of $C$.

\begin{cor}\label{cor3.6}
Let $C\in \TT^{N\times N}$ be an irreducible matrix except for $(+\infty)$, and let $\Lambda$ be the eigenvalue of $C$.
Then, there exists a diagonal matrix $\Gamma\in\mathcal{D}(N)$ such that the matrix  $C':=(-\Gamma)\odot C\odot \Gamma=(X_{i,j}')_{i,j}$ satisfy the equation $\Lambda=\min_{j=1,2,\dots,N}[X_{i,j}']$ for each $i=1,2,\dots,N$.
\end{cor}
\proof Set $C=(X_{i,j})_{i,j}$.
Let $\V{V}=(V_i)_i$ be a tropical eigenvector of $C$.
Then, we have $C\odot \V{V}=\Lambda+\V{V}$, or equivalently, 
$
\min_{j}[X_{i,j}+V_j-V_i]=\Lambda$, ($i=1,2,\dots,N$).
It suffices to define $\Gamma:=\mathrm{diag}[V_1,\dots,V_N]$. $\qed$

\section{Flat-slanted form of weighted digraphs}\label{sec4}

\subsection{Condensation map and Strong-condensation map}\label{sec4.1}

Let $G$ be a weighted digraph.
Introduce two equivalent relations $\approx_w$ and $\approx_s$ on $\mathrm{node}(G)$ by:
\begin{center}
$i\approx_w j\iff$ nodes $i$ and $j$ are connected in $\tp{G}$,\\
$i\approx_sj\iff$ nodes $i$ and $j$ are strongly connected in $\tp{G}$.
\end{center}
Of course,  $i\approx_s j$ implies $i\approx_w j$.

We call the quotient graph $G/\approx_w$ ({\it resp.}\,$G/\approx_s$) the {\it condensation graph of $G$} ({\it resp.}\,{\it strong-condensation graph of $G$}).  
The quotient homomorphism $G\to G/\approx_w$ ({\it resp.}\,$G\to G/\approx_s$) will be called {\it condensation map} ({\it resp.}\,{\it strong-condensation map}).
We note a few of fundamental facts which can be proved immediately:
\begin{itemize}
\item If $s(G)<+\infty$, then $s(G/\approx_w)>s(G)$.
\item If $s(G)<+\infty$, then $(G/\approx_s)_{\leq s(G)}$ is a forest.
\item Any condensation graph contains no simple loop.
\end{itemize}

\begin{lemma}\label{multi}
Let $G$ be a weighted digraph, and $\mu$ be a rational number.
Introduce an equivalence relation $\sim_\mu$ on $\mathrm{node}(G)$ by $i\sim_\mu j\iff $ $i$ and $j$ are connected in $G_{\leq \mu}$.
Then, the homomorphism $\pi:G\to G/\sim_\mu$ is a composition of finitely many condensation maps. \qed
\end{lemma}

\begin{defi}
A composition of finitely many condensation maps is called a {\it multi-condensation map}.
\end{defi}

We denote "$G_1\vartriangleright G_2$" if $G_1\to G_2$ is a condensation map and "$G_1\blacktriangleright G_2$" 
if $G_1\to G_2$ is a multi-condensation map.

\subsection{Similarity translation of weighted graphs}

Let $G$ be a weighted digraph with $N$ nodes.
For a rational number $\nu\in\QQ$ and a node $k\in G$, we define a new weighted digraph $S_{k}(\nu)\cdot G$
\begin{itemize}
\item by replacing $k\xrightarrow{X}i $ with $k\xrightarrow{X+\nu}i $ for $\forall i\in G$ and
\item by replacing $j\xrightarrow{X}k $ with $j\xrightarrow{X-\nu}k $ for $\forall j\in G$.
\end{itemize} 
For any tropical matrix $C$ and its adjacency graph $G(C)$, the following equation holds:
\begin{equation}
S_k(\nu)\cdot G(C)=G((-\gamma)\odot C\odot\gamma),\qquad \gamma=\mathrm{diag}[0,\dots,0,\overset{k}{\hat{\nu}},0,\dots,0]. 
\end{equation}
A {\it similarity translation of weighted digraphs} is a composition of finitely many translations: $G\mapsto S_k(\mu)\cdot G$ ($k\in G$, $\mu\in\QQ$).
We denote $G_1\sim G_2$ if there exists a similarity transformation $G_1\mapsto G_2$.
The relation $\sim$ is an equivalence relation, which we will call the {\it similarity equivalence relation}.

Let $f:G\to H$ be a homomorphism of weighted digraphs.
For $i\in H$ and $\mu\in\QQ$, we define the {\it pull back of the similarity transformation $S_i(\mu)$ by $f$} by
\[
f^\ast S_i(\mu):=S_{j_1}(\mu)\circ\dots\circ S_{j_k}(\mu),\qquad
\mbox{where } f^{-1}(i)=\{j_1,\dots,j_k\}.
\]
(Note that $S_i(\mu)\circ S_j(\nu)=S_j(\nu)\circ S_i(\mu)$). 
For a general similarity transformation $S=S_{i_1}(\mu_1)\circ\dots\circ S_{i_l}(\mu_l)$, 
we set $f^{\ast}S:=f^\ast S_{i_1}(\mu_1)\circ\cdots\circ f^\ast S_{i_l}(\mu_l)$.
There naturally exists an induced homomorphism $f^\ast S\cdot G\to S\cdot H$ whose restriction to the node sets coincides with that of $f:G\to H$.
In other words, 
\begin{equation}\label{diagram}
\mbox{the diagram }  
\begin{array}{ccc}
G  &  & \\
\downarrow & & \\
H & \sim & S\cdot H
\end{array}
\mbox{ implies }
\begin{array}{ccc}
G  & \sim& f^\ast S\cdot G\\
\downarrow & & \downarrow\\
H & \sim & S\cdot H
\end{array}.
\end{equation}
\begin{rem}\label{rem4.1}
If $J\subset G$ is a subgraph such that $f(J)=(\mbox{a single node})$, $J$ is invariant under the similarity transformation $G\sim f^\ast S\cdot G$.
\end{rem}


\begin{lemma}\label{lemma4.2}
Let $A$ be a weighted digraph with at least one arc, $B$ be the condensation or strong-condensation graph of $A$, and $D$ be a digraph which is similarity equivalent to $B$.
Consider the diagram 
$
\begin{array}{ccc}
A  & \sim& \exists C\\ 
\downarrow & & \downarrow\\
B & \sim & D
\end{array}
$
which is derived from (\ref{diagram}). 
In this situation, the following (1--3) hold:
\begin{enumerate}
\renewcommand{\theenumi}{(\arabic{enumi})}
\item $s(A)\leq s(D)$ implies $s(A)\leq s(C)$.
\item If (i) $B$ is the condensation graph of $A$ or (ii)  $B$ is the strong-condensation graph of $A$ and $\tp{A}$ contains at least one loop, then $s(A)\leq s(D)$ implies $s(A)=s(C)$.
\item In addition to (2), if $s(A)<s(D)$, then the leading term of $C$ is expressed as
\[
\tp{C}=\left\{
\begin{array}{cc}
\tp{A}, & \mbox{the case (i)}\\
\mbox{(the disjoint union of strongly connected components of $\tp{A}$)}, & \mbox{the case (ii)}
\end{array}
\right..
\]
\end{enumerate}
\end{lemma}
\proof (1): Let $f:A\to B$ and $g:C\to D$ be the homomorphisms in the diagram.
Consider an arc $c=(i\to j)$.
From Remark \ref{rem4.1}, we have $g(i)=g(j)\Rightarrow w_C(c)=w_A(c)\geq s(A)$.
On the other hand, 
if $g(i)\neq g(j)$, the inequality $w_C(c)\geq w_D(g(i)\to g(j))\geq s(D)\geq s(A)$ holds.
Therefore, $s(A)\leq w_C(c)$, ($\forall c\in C$).
This implies $s(A)\leq s(C)$.
(2): Let
\begin{gather*}
\mathcal{A}:=
\{a=(i\to j)\in A\,\vert\,f(i)=f(j)\mbox{ and }w_A(a)=s(A)\},\\
\mathcal{A}':=
\{c=(i\to j)\in C\,\vert\,g(i)=g(j)\mbox{ and }w_C(c)=s(A)\}.
\end{gather*}
From Remark \ref{rem4.1}, $\mathcal{A}$ and $\mathcal{A}'$ are isomorphic as sets.
When one of the conditions (i) and (ii) is satisfied, we have $\mathcal{A}\neq \emptyset$, which implies $s(A)\geq s(C)$. 
By the first part of this lemma, we conclude $s(A)=s(C)$. 
(3): Because $g(i)\neq g(j)\Rightarrow w_C(i\to j)\geq s(D)>s(A)=s(C)$, we have $\mathcal{A}'=\{c\in C\,\vert\,w_C(c)=s(C)\}$. 
Therefore, $((w_A(i\to j)=s(A))\wedge (f(i)=f(j) ))\iff w_C(i\to j)=s(C)$. 
This concludes (3). 
$\qed$

From Lemma \ref{lemma4.2} (3), 
\begin{equation}\label{eq16}
\mbox{a diagram } 
\begin{array}{ccc}
A  &  & \\
\triangledown & & \\
B & \sim & D
\end{array},\ (s(A)<s(D)) \mbox{ induces the diagram } \begin{array}{ccc}
A  & \sim& C\\
\triangledown & & \triangledown\\
B & \sim & D
\end{array}
,\ (s(A)=s(C)). 
\end{equation}
More generally, we have the following lemma:
\begin{lemma}\label{lemma4.5}
A diagram 
$\begin{array}{ccc}
A  &  & \\
\blacktriangledown & & \\
B & \sim & D
\end{array}$,
($s(A)<s(D)$) induces the diagram 
$\begin{array}{ccc}
A  & \sim& C\\
\blacktriangledown & & \blacktriangledown\\
B & \sim & D
\end{array}$,
($s(A)=s(C)$).
\end{lemma}
\proof Because $A\blacktriangleright B$ is a multi-condensation map, there exists a sequence of condensation maps: $A\vartriangleright A_1\vartriangleright\cdots\vartriangleright A_K\vartriangleright B$.
By definition of condensation maps, we have $s(A)<s(A_1)<\dots<s(A_K)<s(B)$.
Therefore, by using (\ref{eq16}) repeatedly, we can find the diagram 
$\begin{array}{ccc}
A_n  & \sim& \exists C_n\\
\triangledown & & \triangledown\\
A_{n+1} & \sim & C_{n+1}
\end{array}$, 
($s(A_n)=s(C_n)$) for all $n$. $\qed$

\subsection{Definition of flat-slanted graph}

\begin{lemma}\label{lemma4.4}
Let $A$ be a weighted digraph whose leading term contains at least one loop.
There exists a digraph $C$ such that (i) $A\sim C$ and (ii) $\tp{C}$ is a disjoint union of strongly connected components.
\end{lemma}
\proof Set $s=s(A)$.
Let $\pi:A\to (A/\approx_s)=:B$ be the strong-condensation map.
Because $B_{\leq s}$ is a forest (see \S \ref{sec4.1}), we can define the semi-order $\succ$ on $\mathrm{node}(B)$ by $x\succ y\iff$ there exists a path from $x$ to $y$ in $B_{\leq s}$.
Fix an order-preserving bijection $\varphi:\mathrm{node}(B)\to\{1,2,\dots,n\}$ and set $D:=\{S_{x_1}(\epsilon\cdot \varphi(x_1))\circ\dots\circ S_{x_n}(\epsilon\cdot \varphi(x_n))\}\cdot B$ for a rational number $\epsilon>0$.
If $\epsilon>0$ is sufficiently small, we can assume $s(D)>s(A)$.
Let $C$ be the weighted digraph uniquely defined by (\ref{diagram}) as
$
\begin{array}{ccc}
A & \sim &  C \\
\downarrow & & \downarrow\\
B & \sim & D 
\end{array}
$.
From Lemma \ref{lemma4.2} (3), $\tp{C}=C_{\leq s}=$ (the disjoint union of strongly connected components of $\tp{A}$).
$\qed$

\begin{prop}\label{lemma3.7}
Let $G$ be a strongly connected, weighted digraph.
Then, there exists a weighted digraph $H$ such that (i) $G\sim H$ and (ii) $\tp{H}$ is a disjoint union of strongly connected components. 
\end{prop}
\proof If $G=\ast$ (the simple node graph), there is noting to prove.
Assume $G\neq \ast$.
Let $C\neq (+\infty)$ be the adjacency matrix of $G$.
From Corollary \ref{cor3.6}, there exists a diagonal tropical matrix $\Gamma$ such that the matrix $C'=(-\Gamma)\odot C\odot \Gamma=(X_{i,j}')_{i,j}$ satisfies $\Lambda=\min_j[X'_{i,j}]$ for all $i$.
Define $J:=G(C')$.
By the equation above, the digraph $J$ satisfies the following properties : $\bullet$ $G\sim J$, $\bullet$ $s(J)=\Lambda$ and $\bullet$ for each node $i\in J$, there exists at least one arc $i\xrightarrow{\Lambda} j$.
By pigeonhole principle, the digraph $\tp{J}=J_{\leq \Lambda}$ contains at least one loop. 
By Lemma \ref{lemma4.4}, there exists a weighted digraph $H$ such that (i) $J\sim H$ and (ii) $\tp{H}$ is a disjoint union of strongly connected components. $\qed$

\begin{defi}
A weighted digraph whose leading term is a disjoint union of strongly connected components is said to be {\it flat-slanted}.
\end{defi}
The results in this section can be summarized as follows:
\begin{prop}\label{flat}
Any strongly connected weighted digraph is similarity equivalent to a flat-slanted graph.\qed
\end{prop}

Let $H$ be a flat-slanted graph with at least one arc.
Because $s(H)=\Lambda(H)$, we have $H\vartriangleright H'\Rightarrow \Lambda(H)=s(H)<s(H')\leq \Lambda(H')$. 

\begin{prop}\label{prop4.7}
If $D$ is flat slanted, the diagram 
$
\begin{array}{ccc}
A & & \\
\blacktriangledown & & \\
B & \sim & D 
\end{array}
$ implies
$
\begin{array}{ccc}
A & \sim & C\\
\blacktriangledown & &\blacktriangledown \\
B & \sim & D 
\end{array}
$.
\end{prop}
\proof The statement follows from $s(A)<s(B)\leq \Lambda(B)=\Lambda(D)=s(D)$ and Lemma \ref{lemma4.5}. $\qed$

\begin{example}\label{sec4.2.1}

Let 
$C:=
\left[ 
\begin{array}{ccc}
1 & 1 & 0 \\ 
0 & 4 & 2 \\ 
3 & 3 & 2
\end{array} 
\right] \in \TT^{3\times 3}
$
be the tropical matrix, whose adjacency graph is 
$
G=
\left(
\vcenter
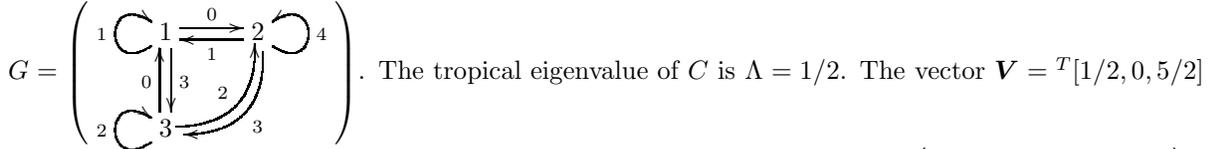
{
\xymatrix{
1 \ar@(ld,lu)^1\ar@<0.5ex>[r]^0\ar@<0.5ex>[d]^3& 2 \ar@(rd,ru)_4\ar@<0.5ex>[l]^1 \ar@<0.5ex>@/^15pt/[dl]^3\\
3 \ar@(ld,lu)^2\ar@<0.5ex>[u]^0\ar@<0.5ex>@/_15pt/[ur]^2
}
}
\right)$.
The tropical eigenvalue of $C$ is $\Lambda=1/2$.
The vector $\V{V}={}^T[1/2,0,5/2]$ is a corresponding eigenvector.
Set $J:=\{S_1(1/2)\circ S_2(0)\circ  S_3(5/2)\}\cdot G=
\left(
\vcenter{
\xymatrix{
1 \ar@(ld,lu)^1\ar@<0.5ex>[r]^{1/2}\ar@<0.5ex>[d]^{1}& 2 \ar@(rd,ru)_4\ar@<0.5ex>[l]^{1/2} \ar@<1ex>@/^15pt/[dl]^{1/2}\\
3 \ar@(ld,lu)^2\ar@<0.5ex>[u]^2\ar@<0ex>@/_15pt/[ur]^(.33){9/2}
}
}
\right)$.
The leading term $\tp{J}=J_{\leq 1/2}$ contains two strongly connected components $x=\{1,2\}$, $y=\{3\}$.
These strongly connected components are connected by one arc $(2\xrightarrow{1/2} 3)$.
Next define $H:=\{S_1(\epsilon)\circ S_2(\epsilon)\}\cdot J$ for sufficiently small $\epsilon>0$.
Then, the graph $H$ is flat-slanted.
For example, if $\epsilon=1/4$, we have
$
H=
\left(
\vcenter
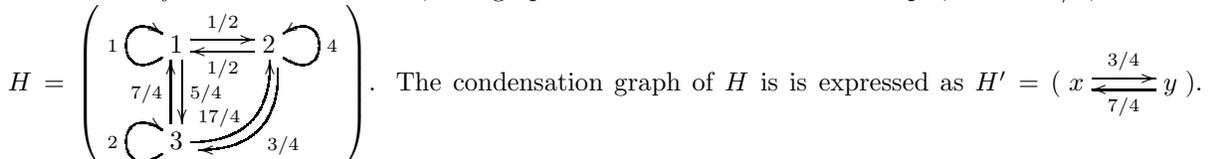
{
\xymatrix{
1 \ar@(ld,lu)^1\ar@<0.5ex>[r]^{1/2}\ar@<0.5ex>[d]^{5/4}& 2 \ar@(rd,ru)_4\ar@<0.5ex>[l]^{1/2} \ar@<1ex>@/^15pt/[dl]^{3/4}\\
3 \ar@(ld,lu)^2\ar@<0.5ex>[u]^{7/4}\ar@<0ex>@/_15pt/[ur]^(.33){17/4}
}
}
\right)$.
The condensation graph of $H$ is is expressed as
$
H'=(
\xymatrix{
x\ar@<0.5ex>[r]^{3/4} & y\ar@<0.5ex>[l]^{7/4}
}
)
$.
Obviously, $\Lambda(G)=\Lambda(J)=\Lambda(H)=1/2<\Lambda(H')=5/4$.
\end{example}

\section{Gently-slanted form of weighted digraphs}\label{sec5}

In this section, we introduce the {\it gently-slanted graph form} of weighed digraphs, obtained by "perturbing" weights of flat-slanted graphs.

\subsection{Topology on the set of weighted digraphs}

To make precise our claims, we need a topology on a set of weighted digraphs.
Let $U$ be an unweighted digraph and $\mathfrak{X}_U$ be the set of weighted digraphs whose underlying graph ({\it i.e.}\,the unweighted graph obtained by forgetting weights) is $U$.
Let $\mathrm{arc}(U)=\{a_1,\dots,a_A\}$ be the arc set of $U$.
We introduce the coarsest topology on $\mathfrak{X}_U$ such that the map
\[
\mathfrak{X}_U\to \QQ^{A};\quad G\mapsto (w_G(a_i))_{i=1,2,\dots,A}
\]
is continuous.
Obviously, the map $\mathfrak{X}_U\times \QQ\to \mathfrak{X}_U;$ $(G,\mu)\mapsto S_i(\mu)\cdot G$ is continuous.

\subsection{Definition of the gently-slanted form}

\begin{lemma}\label{lemma5.1}
Let $A$ be a weighted digraph without simple loops whose leading term is strongly connected.
For any proper and non-empty subset $E\subsetneq \mathrm{node}(A)$, there exist a positive number $\epsilon$ and a continuous map
\[
B_E:(s(A)-\epsilon,s(A))\to \mathfrak{X}_U;\qquad u\mapsto B_E(u)
\]
such that (i) $A\sim B_E(u)$, (ii) $\tp{B_E(u)}$ is a $E$-forest, (iii) $u=s(B_E(u))$ and (iv) $\lim\limits_{u\uparrow s(A)}B_E(u)=A$.
\end{lemma}
\proof 
Let $\delta:=s(A)-u>0$. 
Write $E=\{i_1,\dots,i_e\}$, ($e:=\sharp E$).
For a path $\gamma\in A$ whose starting node is $i_k$, define $l(\gamma):=\mathrm{length}(\gamma)+\frac{k-1}{e}\delta$.
Define the rational number valued function $d:\mathrm{node}(A)\to \frac{1}{e}\cdot \ZZ_{\geq 0}$ by $d(j):=\min_{i\in E}\{l(\gamma)\,\vert\,\gamma \mbox{ is a path from $i$ to $j$ in $\tp{A}$}\}$.
($d$ is well-defined because $\tp{A}$ is strongly connected).
Obviously, $d(i_1)=0$, $d(i_2)=1/e$, $\dots$, $d(i_e)=(e-1)/e$.
Label the nodes of $A$ as $i_1,\dots,i_N$ such that  $1\leq d(i_{e+1})\leq d(i_{e+2})\leq \dots\leq d(i_N)$.
Let
$
X_{m,n}:=
\left\{
\begin{array}{cl}
w_A(i_n\to i_m), & \mbox{if } (i_n\to i_m)\in A\\
+\infty, & \mbox{else}
\end{array}
\right.
$.
Define the rational sequence $\alpha_1,\alpha_2,\dots,\alpha_N$ recursively by
$$
\left\{
\begin{array}{l}
\alpha_1=0,\ 
\alpha_2=\frac{1}{e}\delta,\
\alpha_3=\frac{2}{e}\delta,\dots,
\alpha_e=\frac{e-1}{e}\delta,\\
\alpha_{k+1}:=\min_{l=1}^k[X_{k+1,l}+\alpha_l]-u,\quad (k\geq e)
\end{array}
\right..
$$
If the number $\delta>0$ is sufficiently small, $\alpha_k$ is equal to $d(i_k)\cdot \delta$.
(It is proved by induction).
We will prove that the graph
\begin{equation}\label{bee}
B_E(u):=\{S_{i_1}(\alpha_1)\circ S_{i_2}(\alpha_2)\circ\dots\circ S_{i_N}(\alpha_N)\}\cdot A
\end{equation}
is a desired one.
(i) and (iv) are obvious.
(iii): Write $B=B_E(u)$.
By (\ref{bee}), we have 
$$
w_B(i_m\to i_n)=w_A(i_m\to i_n)+\alpha_m-\alpha_n=X_{n,m}+\alpha_m-\alpha_n\overset{\delta:\mathrm{small}}{=}X_{n,m}+\{d(i_m)-d(i_n)\}\delta.
$$
On the other hand, by definition of $d$, we have ($\dagger$): $X_{n,m}=s(A) \Rightarrow d(i_m)-d(i_n)\geq -1$, and ($\dagger\dagger$): $d(i_n)\geq 1\Rightarrow \exists m<n$ (($X_{n,m}=s(A)$) $\wedge$ ($d(i_m)-d(i_n)=-1$)).
We obtain $u\geq s(B_E(u))$ from ($\dagger$), and $u\leq s(B_E(u))$ from ($\dagger\dagger$).
Hence, (iii) holds.
(ii): Because $(i_m\to i_n)\in \tp{B_E(u)}\iff d(i_m)-d(i_n)=-1$, the digraph $\tp{B_E(u)}$ is decomposed into $e$ connected components.
Each component is an acyclic graph on $i_n$ ($n=1,2,\dots,e$). 
$\qed$

\begin{prop}\label{prop5.3}
Let $H$ be a flat-slanted graph without simple loops and $\tp{H}=h^1\sqcup\dots \sqcup h^n$ be the strongly connected component decomposition. 
Fix a proper subset $E\subsetneq \mathrm{node}(H)$ such that $E\cap h^i\neq \emptyset$ for $i=1,2,\dots,n$.
Then, there exist a positive number $\epsilon>0$ and a continuous map
\[
F_E:(s(H)-\epsilon,s(H))\to\mathfrak{X}_U;\qquad u\mapsto F_R(u)
\]
such that (i) $H\sim F_E(u)$, (ii) $\tp{F_E(u)}$ is a $E$-forest, (iii) $u=s(F_E(u))$ and (iv) $\lim\limits_{u\uparrow s(H)}F_E(u)=H$.
\end{prop}
\proof Set $E^i:=E\cap h^i$.
Let $I:=\{i\,\vert\,E^i\subsetneq h^i\}\neq \emptyset$.
From Lemma \ref{lemma5.1}, there exists a similarity transformation $S^i(u)$, ($i\in I$) such that the digraph  $B_{E^i}(u)=S^i(u)\cdot h^i$ satisfies $\bullet$ $\tp{B_{E^i}(u)}$ is a $E^i$-forest and $\bullet$ $u=s(B_{E^i}(u))$.
Define $F_E(u):=\{S^{i_1}(u)\circ\cdots\circ S^{i_m}(u)\}\cdot H$ for $I=\{i_1,\dots,i_m\}$.
If the difference $s(H)-u>0$ is sufficiently small, we have $\tp{F_E(u)}=\tp{B_{E^{i_1}}(u)}\sqcup\dots\sqcup\tp{B_{E^{i_m}}(u)}\sqcup \ast\sqcup \cdots\sqcup\ast$
($\ast$= a separated node).
Because a disjoint union of $E^i$-forests is a $(E^1\cup\dots\cup E^n)$-forest, we conclude the claim.
$\qed$

\begin{defi}
Let $F$ be a digraph, and $E\subsetneq \mathrm{node}(F)$ be a proper subset.
A weighted digraph $F$ is said to be {\it $E$-gently-slanted} if $\tp{F}$ is a $E$-forest. 
\end{defi} 

Proposition \ref{prop5.3} (i), (ii) implies that {\it $H$ is similarity equivalent to a $E$-gently-slanted graph}.

\begin{example}
Let 
$C:=
\left[ 
\begin{array}{ccc}
+\infty & 1 & 0 \\ 
0 & +\infty & 2 \\ 
3 & 3 & +\infty
\end{array} 
\right]
$
be the tropical matrix, whose adjacency graph is 
$
G=
\left(
\vcenter{
\xymatrix{
1 \ar@<0.5ex>[r]^0\ar@<0.5ex>[d]^3& 2 \ar@<0.5ex>[l]^1 \ar@<0.5ex>@/^15pt/[dl]^3\\
3 \ar@<0.5ex>[u]^0\ar@<0.5ex>@/_15pt/[ur]^2
}
}
\right)
$.
By similar method in Example \ref{sec4.2.1}, we obtain the flat-slanted graph
$
H=
\left(
\vcenter{
\xymatrix{
1 \ar@<0.5ex>[r]^{1/2}\ar@<0.5ex>[d]^{5/4}& 2 \ar@<0.5ex>[l]^{1/2} \ar@<1ex>@/^15pt/[dl]^{3/4}\\
3 \ar@<0.5ex>[u]^{7/4}\ar@<0ex>@/_15pt/[ur]^(.33){17/4}
}
}
\right)$ such that $G\sim H$.
The leading term $\tp{H}=H_{\leq 1/2}$ contains two strongly connected components $x=\{1,2\}$, $y=\{3\}$ and two arcs.
Therefore, there are two possibilities: (i) $E=\{1,3\}$ or (ii) $E=\{2,3\}$.

(i) When $E=\{1,3\}$, $F=S_2(\epsilon)\cdot H$ is a $E$-gently-slanted graph for sufficiently small $\epsilon>0$.
In fact, 
$F=
\left(
\vcenter{
\xymatrix{
\ \ \ \ 1\ \ \ \  \ar@<0.5ex>[r]^{1/2-\epsilon}\ar@<0.5ex>[d]^{5/4}& 2 \ar@<0.5ex>[l]^{1/2+\epsilon} \ar@<1ex>@/^15pt/[dl]^{1+\epsilon}\\
\ \ \ \ 3\ \ \ \  \ar@<0.5ex>[u]^{7/4}\ar@<0ex>@/_15pt/[ur]^(.33){17/4-\epsilon}
}
}
\right)
$ is $E$-gently-slanted.

(ii) When $E=\{2,3\}$, $F=S_1(\epsilon)\cdot H$ is a $E$-gently-slanted graph for sufficiently small $\epsilon>0$.
In fact,
$F=
\left(
\vcenter{
\xymatrix{
\ \ 1\ \  \ar@<0.5ex>[r]^{1/2+\epsilon}\ar@<0.5ex>[d]^{5/4+\epsilon}& 2 \ar@<0.5ex>[l]^{1/2-\epsilon} \ar@<1ex>@/^15pt/[dl]^{3/4}\\
\ \ 3\ \  \ar@<0.5ex>[u]^{7/4-\epsilon}\ar@<0ex>@/_15pt/[ur]^(.33){17/4}
}
}
\right)
$ is $E$-gently-slanted.
\end{example}

\section{The Perron theorem for positive Puiseux matrices}\label{sec6}

\subsection{Real Puiseux field and its subtraction free part}

Let $F=\bigcup_{p=1}^\infty\RR((t^{1/p}))$ be the real Puiseux series field.
The {\it positive part} $F_>$ of $F$ is the subset of positive elements ({\it i.e.}\,elements whose leading term is positive).
We call the set $F_\geq:=F_>\cup\{0\}$ the {\it nonnegative part of} $F$.
$F$ admits the total order $\geq$ defined by $x\geq y \iff x-y\in F_\geq$.
It is known that the ordered field $F$ is a {\it real closed field} (see Section \ref{secaa}).

The {\it subtraction-free part of} $F$ is the subset $\mathfrak{sf}$ of $F_>$ defined by
\[
\mathfrak{sf}:=\{a_{k/q}t^{k/q}+
a_{(k+1)/q}t^{(k+1)/q}+\cdots\,\vert\,k\in \ZZ,\ q\in\ZZ_{>0},\ a_{r}\geq 0\, (\forall r),\ a_{k/q}>0 \}.
\]

A real Puiseux matrix is said to be {\it positive} ({\it resp.\,nonnegative, subtraction-free}) if all its entries are contained in $F_>$ ({\it resp.}\,$F_\geq$, $\mathfrak{sf}$).
An element of $F^N$ ({\it resp.}\,$F_>^N$, $F_\geq^N$) is called a {\it real Puiseux vector} ({\it resp.\,a positive, nonnegative Puiseux vector}).

The following is an analog of the Perron theorem over the real Puiseux field.

\begin{thm}[The Perron theorem over $F$]\label{main}
A positive Puiseux matrix $Y$ satisfies the following properties (i--iv): (i) Its spectral radius $\rho(Y)$ is a simple and positive eigenvalue. (ii) A corresponding eigenvector with $\rho(Y)$ can be chosen to be entry-wise positive. (iii) Any nonnegative eigenvector corresponds with $\rho(Y)$. (iv) No eigenvalue $\lambda$, except for $\rho(Y)$, satisfies $\zet{\lambda}=\rho(Y)$.
\end{thm}
We will refer to $\rho(Y)$ as the {\it Perron root} of $Y$.
A {\it Perron vector} is a positive Puiseux eigenvector corresponding to the Perron root.


In the sequel, we give a constructive proof of Theorem \ref{main} for some generic class of Puiseux matrices (see \S \ref{sec.generic}) by demonstrating the method of calculating Perron roots and Perron vectors of positive Puiseux matrices.

\begin{rem}
Any positive Puiseux matrix $Y\in F_>^{N\times N}$ can be decomposed as $Y=kY'$, where $k\in F_>$ and $Y'\in \mathfrak{sf}^{N\times N}$.
Hereafter, we can restrict ourselves on the subtraction-free matrices only.
\end{rem}

\subsection{Perron-Frobenius data}

\begin{defi}
An EQ $\mathcal{X}=(Y;\lambda;\{\V{x}\},\{\V{y}\})$ of depth $\Lambda$ is said to have the {\it Perron-Frobenius property} if it satisfies the following (i--iv):
\begin{enumerate}
\renewcommand{\theenumi}{(\roman{enumi})}
\item $\lambda\equiv \rho(Y)\mod{I(v(Y-\lambda\cdot\mathrm{id})+\Lambda)}$,
\item The quasi basis $\{\V{x}\}$ is expressed as $\{\V{x}\}=\{\V{x}^1\,t^\emptyset\V{x}^2,\dots,t^\emptyset\V{x}^N\}$, and the quasi-basis $\{\V{y}\}$ is expressed as $\{\V{y}\}=\{\V{y}^1,t^\emptyset\V{y}^2,\dots,t^\emptyset\V{y}^N\}$, where $\V{x}^1\in F_{\geq }^N$, $\V{y}^1\in F_{\geq}^N$ and ${}^T\tp{(\V{y}^1)}\cdot \tp{(\V{x}^1)}=1$.
\item Any nonnegative $\Lambda$-th approximate eigenvector of $Y$ corresponds with $\rho(Y)$.
\item No $\Lambda$-th approximate eigenvalue $\lambda$ of $Y$, except for $\rho(Y)$, satisfies  $\zet{\lambda}\equiv\rho(Y)\mod{I(v(Y-\lambda\cdot\mathrm{id})+\Lambda)}$. 
\end{enumerate}
We call the vector $\V{x}^1$ a {\it $\Lambda$-th approximate Perron vector}.
\end{defi}

\begin{rem}\label{depth0}
If $\mathcal{X}=(Y;\lambda;\{\V{x}\};\{\V{y}\})$ is an EQ with Perron-Frobenius property of depth $0$, then $\lambda$ is the Perron root of the real matrix $\tp{Y}$.
The quasi-basis $\{\V{x}\}$ consists of one element $\{\V{x}\}=\{\V{x}^1\}$, where $\V{x}^1\in \RR^N$ is a (usual) Perron vector of $\tp{Y}$.
\end{rem}

For a real Puiseux matrix $X\in F^{M\times N}$ and subsets $\alpha\subset\{1,2,\dots,M\}$,  $\beta\subset\{1,2,\dots,N\}$, we denote $X[\alpha,\beta]$ the submatrix of $X$ whose rows are indexed by $\alpha$ and whose columns are indexed by $\beta$.
If $M=N$ and $\alpha=\beta$, write $X[\alpha]:=X[\alpha,\alpha]$.

\begin{defi}
Let $Y\in\mathfrak{sf}^{N\times N}$ be a subtraction-free matrix and $\alpha=(\alpha_1,\dots,\alpha_m)$ be a partition of $\{1,2,\dots,N\}$ ({\it i.e.}\,a collection of subsets $\alpha_i\subset \{1,2,\dots,N\}$ such that  $\alpha_i\cap\alpha_j=\emptyset$ and $\bigcup_{i=1}^m\alpha_i=\{1,2,\dots,N\}$).
A tuple 
\begin{equation}\label{data}
\mathfrak{A}:=(Y;\mathcal{X}(\alpha_1),\dots,\mathcal{X}(\alpha_m)),\quad\mbox{where}\ 
\mathcal{X}(\alpha_i)=(Y[\alpha_i];\lambda(\alpha_i);\{\V{x}(\alpha_i)\},\{\V{y}(\alpha_i)\})
\end{equation}
is said to be a {\it quasi Perron-Frobenius data} (or shortly, a {\it quasi PF-data}) {\it associated with} $(Y,\alpha)$ if each $\mathcal{X}(\alpha_i)$ is an EQ with the Perron-Frobenius property.
We define the {\it depth of the quasi PF-data} as $\mathrm{depth}(\mathfrak{A}):=\min_i[\mathrm{depth}(\mathcal{X}(\alpha_i))]$.
\end{defi}

\begin{defi}
A quasi PF-data (\ref{data}) is said to be {\it singular} if $\lambda(\alpha_1)=\lambda(\alpha_2)=\dots=\lambda(\alpha_m)$.
\end{defi}

The {\it adjacency graph of a quasi PF-data $\mathfrak{A}$} is the weighted digraph $G(\mathfrak{A})=(\mathrm{node}(G(\mathfrak{A})),\mathrm{arc}(G(\mathfrak{A})))$ which is defined as follows:
\begin{gather*}
\mathrm{node}(G(\mathfrak{A})):=\{\alpha_1,\dots,\alpha_m\},\qquad
(\alpha_i\xrightarrow{W}\alpha_j)\in\mathrm{arc}(G(\mathfrak{A}))\iff W=v\left(Y[\alpha_j,\alpha_i]\right).
\end{gather*}


\begin{defi}
A quasi PF-data $\mathfrak{A}$ (\ref{data}) associated with $(Y,\alpha)$ is said to be a {\it Perron-Frobenius data} (or shortly, a {\it PF-data}) if the Puiseux matrix $Y$ is decomposed as 
\begin{gather}\label{eq16a}
Y\equiv Y[\alpha_1]\oplus\cdots\oplus Y[\alpha_m]+ Rt^{s(G(\mathfrak{A}))}\mod{I(s(G(\mathfrak{A})))},
\end{gather}
where $R\in \RR^{N\times N}$.
\end{defi}

\begin{rem}
From (\ref{eq16a}), we have $G(\mathfrak{A})=G/\sim_{s(G(\mathfrak{A}))}$ if $\mathfrak{A}$ is a PF-data.
See Lemma \ref{multi}.
\end{rem}


\subsection{Outline of the construction}

In the sequel, we will construct a PF-data associated with a subtraction-free matrix of arbitrarily large depth. 
The construction consists of the following three processes:
\begin{enumerate}
\item {\bf Process A}:
   \begin{itemize}
   \item Input: A pair $(Y,\mathfrak{A})$, where $Y$ is a subtraction-free matrix, $\alpha$ is a partition of $\{1,\dots,N\}$ and $\mathfrak{A}$ is a PF-data associated with $(Y,\alpha)$.
   \item Output: A pair $(Z,\mathfrak{B})$, where $Z$ is a subtraction-free matrix and $\mathfrak{B}$ is a PF-data associated with $(Z,\alpha)$ such that: (i) $Y\sim Z$, (ii) $G(\mathfrak{A})\sim G(\mathfrak{B})$, (iii) $G(\mathfrak{B})$ is flat-slanted and (iv) $\mathrm{depth}(\mathfrak{A})=\mathrm{depth}(\mathfrak{B})$.
   \end{itemize}
\item {\bf Process B}:
   \begin{itemize}
   \item Input: A pair $(Y,\mathfrak{A})$, where $Y$ is a subtraction-free matrix, $\alpha$ is a partition of $\{1,\dots,N\}$ and $\mathfrak{A}$ is a PF-data associated with $(Y,\alpha)$ such that (i) $\mathfrak{A}$ is non-singular and (ii) $\tp{G(\mathfrak{A})}$ is strongly-connected.
   \item Output: A pair $(Y,\mathfrak{B})$, where $\mathfrak{B}$ is a PF-data associated with $(Z,\beta)$, $\beta$ is a courser partition than $\alpha$ such that $\mathrm{depth}(\mathfrak{B})>\mathrm{depth}(\mathfrak{A})$.
   \end{itemize}
\item {\bf Process C}:
   \begin{itemize}
   \item Input: A pair $(Y,\mathfrak{A})$, where $Y$ is a subtraction-free matrix, $\alpha$ is a partition of $\{1,\dots,N\}$ and $\mathfrak{A}$ is a PF-data associated with $(Y,\alpha)$ such that (i) $\mathfrak{A}$ is singular, (ii) $\mathfrak{A}$ is of depth $0$ and (iii) $\tp{G(\mathfrak{A})}$ is strongly-connected.
   \item Output: A pair $(Y,\mathfrak{B})$, where $\mathfrak{B}$ is a PF-data associated with $(Z,\{1,\dots,N\})$ such that $\mathrm{depth}(\mathfrak{B})>\mathrm{depth}(\mathfrak{A})$.
   \end{itemize}
\end{enumerate}

The rough flow of the construction is as follows.
From Proposition \ref{flat}, it is sufficient to consider a subtraction-free matrix $Y$ whose adjacent matrix is flat-slanted.

\begin{algo}\label{algo}
Let $Y$ be a subtraction-free matrix whose adjacent matrix is flat-slanted.
\begin{enumerate}
\item[0-a.] Since the adjacent matrix of $Y$ is flat-slanted, $Y$ is decomposed as $Y=Y[\alpha_1]\oplus\dots \oplus Y[\alpha_m]+o(t^0)$, where each $Y[\alpha_i]$ is irreducible.
Define $\alpha:=(\alpha_1,\dots,\alpha_m)$ and $\mathcal{X}(\alpha_i):=(Y[\alpha_i];\lambda(\alpha_i),\{\V{x}\},\{\V{y}\})$, where $\lambda(\alpha_i)$ is the (usual) Perron root of $Y[\alpha_i]$, $\V{x}$ ({\it resp.}\,$\V{y}$) is a Perron vector of $Y[\alpha_i]$ ({\it resp.}\,${}^TY[\alpha_i]$) such that ${}^T\V{y}\cdot \V{x}=1$.
\item[0-b.] Define $\mathfrak{A}:=(Y;\mathcal{X}(\alpha_1),\dots,\mathcal{X}(\alpha_m))$, that is a PF-data associated with $(Y,\alpha)$ of depth $0$.
\item Operate {\bf Process A} on the pair $(Y,\mathfrak{A})$ and obtain a new pair $(Z,\mathfrak{B})$.
Then, update $(Y,\mathfrak{A})\leftarrow (Z,\mathfrak{B})$.
\item Since $G(\mathfrak{A})$ is flat-slanted, $Y$ is re-decomposed as 
\begin{equation}\label{eqre}
Y=Y[\alpha_1]\oplus\dots\oplus Y[\alpha_n]+Rt^{s(G(\mathfrak{A}))}+o(t^{s(G(\mathfrak{A}))})=Y[\beta_1]\oplus\cdots \oplus Y[\beta_l]+o(t^{s(G(\mathfrak{A}))}),
\end{equation}
where $\beta=(\beta_1,\dots,\beta_l)$ is a courser partition than $\alpha$. 
Define $\mathfrak{A}(\beta_i):=(Y[\beta_i];\{\mathcal{X}(\alpha_{k})\,\vert\,\alpha_k\in\beta_i\})$.
\item For each $\beta_i$, operate {\bf Process B} on $(Y[\beta_i],\mathfrak{A}(\beta_i))$ if $\mathfrak{A}(\beta_i)$ is non-singular, or operate {\bf Process C} if $\mathfrak{A}(\beta_i)$ is singular of depth $0$.
Then, obtain a new PF-data $\mathfrak{B}(\beta_i)$ and update $\mathfrak{A}(\beta_i)\leftarrow \mathfrak{B}(\beta_i)$.
\item Update $\mathfrak{A}\leftarrow(Y;\mathcal{X}(\beta_1),\dots,\mathcal{X}(\beta_n))$.
\item Repeat 1--4.
\end{enumerate}
\end{algo}

\subsection{The condition of genericness}\label{sec.generic}

Before proceeding to the proof of Theorem \ref{main}, we note the fact that Algorithm \ref{algo} fails if the decomposition (\ref{eqre}) in the step 2 contains some $\beta_i$ such that $\mathfrak{A}(\beta_i)$ is singular and of depth greater than $0$.
Unfortunately, we would not deal with the case here because singular PF-datum of higher depth possess quite complicated behavior.
In other words, we assume the following condition:

{\bf Condition of genericness}: At each step of Algorithm \ref{algo}, no singular PF-data of depth greater than $0$ is derived.

Note that any linearly perturbed matrix $Y_0+tY_1$ satisfies this condition.

\subsection{Process A}

Let $Y=(x_{i,j})\in \mathfrak{sf}^{N\times N}$ be a square subtraction-free matrix. 
The {\it valuation matrix} of $Y$ is the tropical matrix defined by $\vartheta(Y):=(\vartheta(x_{i,j}))\in\TT^{N\times N}$.
If there is no chance of confusion, we use the abbreviation $G(Y)=G(\vartheta(Y))$, where $G(\vartheta(Y))$ is the adjacency graph of $\vartheta(Y)$.

Assume $G(Y)\sim H$ for some digraph $H$.
By definition, there exists some $r_1,\dots,r_N\in \QQ$ such that $H=\{S_1(r_1)\circ\cdots\circ S_N(r_N)\}\cdot G(Y)$.
We define a subtraction-free matrix $Y^{H/G(Y)}:=\delta^{-1}Y\delta$, where $\delta=\mathrm{diag}(t^{r_1},\dots,t^{r_N})$.
It is easily proved that $G(Y^{H/G(Y)})=H$.

\begin{lemma}\label{lemma6.7}
Let $G=G(Y)$ be the adjacency matrix of $\vartheta(Y)$ and $G(\mathfrak{A})$ be the adjacency graph of a PF-data $\mathfrak{A}$ associated with $Y$.
Then, there exists a multi-condensation map $G\blacktriangleright G(\mathfrak{A})$ of weighted digraphs.
\end{lemma}
\proof Define the map $\pi:\mathrm{node}(G)\to \mathrm{node}(G(\mathfrak{A}))$ by $\pi(x)=\alpha_i\iff x\in \alpha_i$.
From (\ref{eq16a}), the map $\pi$ induces the quotient map $G\to G(\mathfrak{A})=G/\sim_{s(G{(\mathfrak{A})})}$ (see Lemma \ref{multi}).
Therefore, this is a multi-condensation map. \qed

Consider a PF-data $\mathfrak{A}=(Y;\mathcal{X}(\alpha_1),\dots,\mathcal{X}(\alpha_m))$.
Let $\mathcal{H}$ be a weighted digraph with $G(\mathfrak{A})\sim \mathcal{H}$.
Therefore, there exists a digraph $H$ obtained from the diagram (\ref{diagram}):
$
\begin{array}{ccc}
G & \sim & H\\
\downarrow & &\downarrow \\
G(\mathfrak{A})& \sim & \mathcal{H}
\end{array}
$.
Define $\mathfrak{A}^{\mathcal{H}/G(\mathfrak{A})}:=(Y^{H/G};\mathcal{X}(\alpha_1),\dots,\mathcal{X}(\alpha_m))$.
Because $Y^{H/G}[\alpha_i]=Y[\alpha_i]$ (see Remark \ref{rem4.1}), the tuple $\mathfrak{A}^{\mathcal{H}/G(\mathfrak{A})}$ satisfies the conditions in (\ref{data}).

\begin{lemma}
If $s(G)<s(\mathcal{H})$, then $\mathfrak{A}^{\mathcal{H}/G(\mathfrak{A})}$ is a PF-data whose adjacency matrix is $\mathcal{H}$.
\end{lemma}
\proof By Lemmas \ref{lemma4.5} and \ref{lemma6.7}, the map $H\to \mathcal{H}$ is a multi-condensation map.
This implies $Y^{H/G}- Y[\alpha_1]\oplus\dots\oplus Y[\alpha_m]=o(t^{s(\mathcal{H})})$.
$\qed$

\begin{cor}
For any subtraction-free square matrix $Y$ and a PF-data $\mathfrak{A}$ associated with $Y$, there exists a new subtraction free matrix $Z$ and a PF-data $\mathfrak{B}$ associated with $Z$ such that: (i) $Y\sim Z$, (ii) $G(\mathfrak{A})\sim G(\mathfrak{B})$, (iii) $G(\mathfrak{B})$ is flat-slanted and (iv) $\Lambda(\mathfrak{A})=\Lambda(\mathfrak{B})$.
\end{cor}
\proof It is enough to define $\mathfrak{B}:=\mathfrak{A}^{\mathcal{H}/G(\mathfrak{A})}$, where $\mathcal{H}$ is flat-slanted. \qed

\subsection{Process B}

Let $\mathfrak{A}=(Y;\mathcal{X}(\alpha_1),\dots,\mathcal{X}(\alpha_m))$ be a PF-data where $\tp{G(\mathfrak{A})}$ is strongly connected.
Without loss of generality, we can assume  $\lambda(\alpha_1)=\cdots=\lambda(\alpha_k)>\lambda(\alpha_{k+1})\geq\cdots\geq \lambda(\alpha_m)$.

Let $E:=\{\alpha_1,\dots,\alpha_k\}$.
From Proposition \ref{prop5.3}, there exists a $E$-gently-slanted digraph $\mathcal{H}$ with $G(\mathfrak{A})\sim\mathcal{H}$.
Define $(Z,\mathfrak{B}):=(Y^{\mathcal{H}/G(\mathfrak{A})},\mathfrak{A}^{\mathcal{H}/G(\mathfrak{A})})$.
Remember that the weights of the gently-slanted graph depend on a small parameter $\epsilon$ and $s(G(\mathfrak{B}))=s(G(\mathfrak{A}))-\epsilon$, $\lim\limits_{\epsilon\downarrow 0}\mathcal{H}=G(\mathfrak{A})$.

By definition of $E$-forest, $k$ is the number of connected components of $\tp{G(\mathfrak{B})}$.
The matrix $Z$ is decomposed as $Z=Z_1\oplus\cdots\oplus Z_k+o(t^{s(G(\mathfrak{B}))}) $, where each submatrix $Z_i$ is associated with the connected component $\mathcal{C}_i$ of $\tp{G(\mathfrak{B})}$ with $\alpha_i\in\mathcal{C}_i$.

The submatrix $Z_i$ is decomposed as
\begin{align*}
Z_i&=
\left(
\begin{array}{c|c|c|c}
Y[\alpha_i]&  & &\\\hline
 & Y[\alpha_{z_1}]  & & \\\hline
  &   &\ddots & \\\hline 
  &   & & Y[\alpha_{z_q}]
\end{array}
\right)+
R_it^{s(G(\mathfrak{A}))-\epsilon}+o(t^{s(G(\mathfrak{A}))}),
\end{align*}
where $\mathcal{C}_i=\{\alpha_1,\alpha_{z_1}, \cdots, \alpha_{z_q}\}$ and $R_i\neq O$ is a non-negative real matrix.

\begin{lemma}
Let $R_i[\alpha,\beta]$ be the submatrix of $R_i$ associated with index sets $\alpha$, $\beta$.
Then, $R_i[\alpha_i,\alpha_{z_u}]=O$ for all $u$.
\end{lemma}
\proof This is a consequence of the fact that $\tp{G(\mathfrak{B})}$ contains no arrow whose end node is $\alpha_i$.\qed

In other words, $R_i$ is expressed as 
\[
R_i=\left(
\begin{array}{c|c}
 O & O \\\hline
 P_i & \ast
\end{array}
\right).
\]

We denote $W_i:=Y[\alpha_{z_1}]\oplus \cdots\oplus Y[\alpha_{z_q}]$.
Since all the eigenvalues of $Y[\alpha_{z_i}]$ are smaller than $\lambda(\alpha_i)$, the matrix $(W-\lambda(\alpha_i)\cdot \mathrm{Id})$ is regular.
Define 
$$
\V{p}_i:=\V{x}^1(\alpha_i)\oplus (\lambda(\alpha_i)-W)^{-1}P_i\V{x}^1(\alpha_i)t^{s(G(\mathfrak{A})-\epsilon)},
$$
where $\V{x}^1(\alpha_i)$ be the approximate Perron vector of $Y[\alpha_i]$.
It is soon checked that the vector $\V{p}_i$ satisfies
\begin{equation}\label{eq26a}
Z_i\V{p}_i=\lambda(\alpha_i)\V{p}_i+o(t^{s(G(\mathfrak{A}))}).
\end{equation}
From general arguments about real closed fields, we can prove the matrix $(\lambda(\alpha_i)-W)^{-1}$ and the vector $\V{p}_i$ are nonnegative.
(See Lemma \ref{lemmaa4}).

Similarly, we can prove 
\begin{equation}\label{eq27}
{}^TZ_i\V{q}_i=\lambda(\alpha_i)\V{q}_i+o(t^{s(G(\mathfrak{A}))}),
\end{equation}
where $\V{q}_i=\V{y}^1(\alpha_i)\oplus \V{0}$.

Define the partition $\beta=(\beta_1,\dots,\beta_k)$ by $\beta_i:=\alpha_i\cup \alpha_{z_1}\cup\dots\cup \alpha_{z_q}$.
Among $\sharp \beta_i$ approximate eigenvectors of $Y[\alpha_1]\oplus W$, pick up $(\sharp\beta_i-1)$ vectors $\V{x}^2(\beta_i),\V{x}^3(\beta_i),\dots,\V{x}^{\sharp \beta_i}(\beta_i)$ arbitrarily, and similarly, among $\sharp \beta_i$ approximate eigenvectors of ${}^TY[\alpha_1]\oplus {}^TW$, pick up $(\sharp\beta_i-1)$ vectors $\V{y}^2(\beta_i),\V{y}^3(\beta_i),\dots,\V{y}^{\sharp \beta_i}(\beta_i)$ arbitrarily.

It is directly checked that the following quintuple is an EQ of depth $\Lambda+\epsilon$:
\[
\mathcal{X}(\beta_i):=(Z_i;\lambda(\alpha_i),\{\V{x}(\beta_i)\},\{\V{y}(\beta_i)\}),
\]
where $\{\V{x}(\beta_i)\}:=\{\V{p}_i,t^{\epsilon+\emptyset}\V{x}^2(\beta_i),\dots,t^{\epsilon+\emptyset}\V{x}^{\sharp\beta_i}(\beta_i)\}$ and $\{\V{y}(\beta_i)\}:=\{\V{q}_i,t^{\epsilon+\emptyset}\V{y}^2(\beta_i),\dots,t^{\epsilon+\emptyset}\V{y}^{\sharp\beta_i}(\beta_i)\}$.

From (\ref{eq26a}), the tuple $\mathfrak{B}:=(Z;\mathcal{X}(\beta_1),\dots,\mathcal{X}(\beta_k))$ is a PF-data associated with the pair $(Z,\beta)$ of depth $\mathrm{depth}(\mathfrak{A})+\epsilon$.

\subsection{Process C}\label{sec6.5}

Let $\mathfrak{A}=(Y;\mathcal{X}(\alpha_1),\dots,\mathcal{X}(\alpha_m))$ be a PF-data of depth $0$, where $\tp{G(\mathfrak{A})}$ is strongly connected and $\lambda(\alpha_1)=\dots=\lambda(\alpha_m)(=:\lambda)$.

Since $\mathcal{X}(\alpha_i)$ is of depth $0$, we have $\mathcal{X}(\alpha_i)=(Y[\alpha_i],\lambda;\{\V{x}^i\},\{\V{y}^i\})$, where $\V{x}^i,\V{y}^i\in \RR^{n_i}$ are entry-wise positive vectors (see Remark \ref{depth0}). 
Define the piece-wise nonnegative vectors $\V{X}^i$ and $\V{Y}^i$ by
\begin{gather*}
\V{X}^i=
\V{0}\oplus\dots\oplus \V{0}\oplus \V{x}^i\oplus \V{0}\oplus\dots\oplus \V{0},\qquad  (1\leq i\leq m),\\
\V{Y}^i=
\V{0}\oplus\dots\oplus \V{0}\oplus{\V{y}}^i\oplus \V{0}\oplus\dots\oplus \V{0},\qquad  (1\leq i\leq m).
\end{gather*}
Then, the quadruple
\[
\mathcal{X}:=(Y;\lambda;\{\V{X}\},\{\V{Y}\}).
\]
is an EQ of depth $0$.

\begin{lemma}
The principal matrices $\Delta(\mathcal{X})$ and $\Omega(\mathcal{X})$ are of the form:
\begin{gather*}
\Delta(\mathcal{X})=(\mbox{an irreducible matrix whose entries except for diagonal elements are non-negative}),\\
\Omega(\mathcal{X})=\mathrm{Id}.
\end{gather*}
\end{lemma}
\proof 
Without loss of generality, the valuation of $Y$ can be assumed to be $0$.
The matrix $Y$ is expanded as
\[
Y=Y_0+t^{s(G(\mathfrak{A}))}Y'+\cdots.
\]
For $i\neq j$, the $(i,j)$-entry of $\Delta(\mathcal{X})$ is calculated as ${}^T\V{y}^j Y'[\alpha_j,\alpha_i]\V{x}^i$, which is nonnegative.
Moreover, it follows that
\[
{}^T\V{y}^j Y'[\alpha_j,\alpha_i]\V{x}^i\neq 0\iff Y'[\alpha_j,\alpha_i]\neq O\iff \tp{G(\mathfrak{A})}\mbox{ contains } (i\to j).
\]
This relation implies that $\Delta(\mathcal{X})$ is irreducible if and only if $\tp{G(\mathfrak{A})}$ is strongly connected.
This implies the first equation.
The second equation follows from the definition of EQs.
$\qed$

By the classical Perron-Frobenius theorem (Theorem \ref{frob}), the matrix equations $\Delta(\mathcal{X})\V{c}=\mu \Omega(\mathcal{X})\V{c}$ and ${}^T\Delta(\mathcal{X})\V{d}=\mu {}^T\Omega(\mathcal{X})\V{d}$ have the Perron root $\mu^\ast>0$ and Perron vectors with ${}^T\V{d}\cdot \V{c}=1$.
By Corollary \ref{cor2.13}, we obtain the new EQ
\[
\mathcal{X}_+=(Y;\lambda+\mu^\ast t^{s(G(\mathfrak{A}))};\{\V{x}^1,t^\emptyset\V{x}^2,\dots,t^\emptyset\V{x}^g\},\{\V{y}^1,t^\emptyset\V{y}^2,\dots,t^\emptyset\V{y}^g\}),
\]
where
\begin{gather*}
\tp{(\V{x}^1)}=(\V{x}^1,\dots,\V{x}^m)\cdot \V{c},\quad
\tp{(\V{y}^1)}=(\V{y}^1,\dots,\V{y}^m)\cdot \V{d}.
\end{gather*}
It is soon proved that $\mathcal{X}$ has the Perron-Frobenius property, and therefore, the pair $(Y,\mathcal{X})$ is a PF-data of depth $s(G(\mathfrak{A}))>0$.

\section{Examples}\label{sec7}

\subsection{Example I}

Let us consider the positive Puiseux matrix 
$Y=\left(
\begin{array}{ccc}
a & bt & ct^2\\
dt^2 & et^3 & f\\
g & h & kt
\end{array}
\right)$,
where $a,b,c,d,e,f,g,h,k>0$, whose valuation matrix is
$
\left[
\begin{array}{ccc}
0 & 1 & 2 \\
2 & 3 & 0 \\
0 & 0 & 1
\end{array}
\right]
$.
Its adjacency matrix is
$
G=
\left(
\vcenter{
\xymatrix{
1 \ar@(ld,lu)^0\ar@<0.5ex>[r]^2\ar@<0.5ex>[d]^0& 2 \ar@<0.5ex>[l]^1\ar@(rd,ru)_3 \ar@<0.5ex>@/^15pt/[dl]^0\\
3 \ar@<0.5ex>[u]^2\ar@<0.5ex>@/_15pt/[ur]^0\ar@(ld,lu)^1
}
}
\right)$.
Define $H:=S_1(1/2)\cdot G=
\left(
\vcenter{
\xymatrix{
1 \ar@(ld,lu)^{0}\ar@<0.5ex>[r]^{5/2}\ar@<0.5ex>[d]^{1/2}& 2 \ar@<0.5ex>[l]^{1/2}\ar@(rd,ru)_3 \ar@<0.5ex>@/^15pt/[dl]^0\\
3 \ar@<0.5ex>[u]^{3/2}\ar@<0.5ex>@/_15pt/[ur]^0\ar@(ld,lu)^1
}
}
\right)$.
Hence, the digraph $H$ is flat-slanted, whose leading term is a disjoint union of strongly connected components $\alpha_1=\{1\}$ and $\alpha_2=\{2,3\}$.

Let $Z:=Y^{H/G}=\left(
\begin{array}{ccc}
a & bt^{1/2} & ct^{3/2}\\
dt^{5/2} & et^3 & f\\
gt^{1/2} & h & kt
\end{array}
\right)$.
We have $Z[\alpha_1]=(a)$ and $Z[\alpha_2]=\left(
\begin{array}{cc}
et^3 & f\\
h & kt
\end{array}
\right)$.
From the classical Perron-Frobenius theorem, we can calculate EQs $\mathcal{X}[\alpha_1]$, $\mathcal{X}[\alpha_2]$ as $\mathcal{X}[\alpha_1]=(a;a;\{1\},\{1\};0)$ and 
$
\mathcal{X}[\alpha_2]=
\left(
Z[\alpha_2];\sqrt{fh};
\left\{
\left(
\begin{array}{c}
\sqrt{f} \\ 
\sqrt{h}
\end{array} 
\right)
\right\},
\left\{
\dfrac{1}{2\sqrt{fh}}
\left(
\begin{array}{c}
\sqrt{h} \\ 
\sqrt{f}
\end{array} 
\right)
\right\};
0
\right)
$.
Then, $\mathfrak{A}=(Y;\mathcal{X}[\alpha_1],\mathcal{X}[\alpha_2])$ is a PF-data of depth $0$.

The behavior of a PF-data of greater depth depends on the magnitudes of $a$ and $\sqrt{fh}$

\subsubsection{The case if $a=\sqrt{fh}$}
Since ${}^T\V{y}^1[\alpha_2]Y[\alpha_2,\alpha_1]\V{x}^1[\alpha_1]=\frac{g}{2\sqrt{h}}t^{1/2}+\frac{d}{2\sqrt{f}}t^{5/2}$ and ${}^T\V{y}^1[\alpha_1]Y[\alpha_1,\alpha_2]\V{x}^1[\alpha_2]=b\sqrt{f}t^{1/2}+c\sqrt{h}t^{3/2}$, we have
$G(\mathfrak{A})=\left(
\vcenter{
\xymatrix{
\alpha_1 \ar@<0.5ex>[r]^{1/2}& \alpha_2 \ar@<0.5ex>[l]^{1/2}
}
}
\right)$.
Because $G(\mathfrak{A})$ is already flat-slanted, we only need to consider one strongly connected component $\beta=\{\alpha_1,\alpha_2\}$.

Let $\varphi:\mathrm{Im}(Y_0-a\cdot \mathrm{Id})\to \CC^3$ be the linear map defined by $\V{e}_2-\sqrt{\frac{h}{f}}\V{e}_3\mapsto \frac{1}{f}\V{e}_3$.
Define 
\[
\mathcal{X}=\mathcal{X}(\beta)=
\left(
Z[\alpha_1]\oplus Z[\alpha_2];a;
\left\{
(1)\oplus \V{0},
(0)\oplus 
\left(
\begin{array}{c}
\sqrt{f}\\
\sqrt{h}
\end{array}
\right)
\right\}
,
\left\{
(1)\oplus \V{0},
(0)\oplus \frac{1}{2\sqrt{fh}}
\left(
\begin{array}{c}
\sqrt{h}\\
\sqrt{f}
\end{array}
\right)
\right\}
\right).
\]
Then we have 
$
P(\mathcal{X})=
\left(
\begin{array}{cc}
0 & b\sqrt{f} \\ 
\frac{g}{2\sqrt{h}} & 0
\end{array} 
\right)
$ and $Q(\mathcal{X})=\mathrm{Id}$.
By the classical Perron-Frobenius theorem, the matrix equation $P(\mathcal{X})\V{c}=\mu Q(\mathcal{X})\V{c}$ has the Perron root $\mu=\sqrt{\dfrac{gb}{2}\sqrt{\dfrac{f}{h}}}$.
It is soon checked that the vector $\V{c}={}^T(\sqrt{2b\sqrt{fh}},\sqrt{g})$ is a corresponding eigenvector.
On the other hand, the transposed equation ${}^TP(\mathcal{X})\V{d}=\mu {}^TQ(\mathcal{X})\V{d}$ has a eigenvector $\V{d}=\dfrac{1}{2\sqrt{2bg\sqrt{fh}}}{}^T(\sqrt{g},\sqrt{2b\sqrt{fh}})$.
Let 
\begin{gather*}
\V{x}^1_0:=
\left(
\begin{array}{cc}
1 & 0 \\ 
0 & \sqrt{f} \\ 
0 & \sqrt{h}
\end{array} 
\right)\V{c}=
\left(
\begin{array}{c}
\sqrt{2b\sqrt{fh}} \\ 
\sqrt{fg} \\ 
\sqrt{hg}
\end{array} 
\right),\quad 
\V{x}^1_{1/2}:=-\varphi((Z_{1/2}-\mu)\V{x}^1_0)=g\sqrt{\dfrac{b}{2\sqrt{fh}}}\,\V{e}_3,\quad
\V{x}^2_0:=\V{e}_1\\
\V{y}^1_0:=
\left(
\begin{array}{cc}
1 & 0 \\ 
0 & \frac{1}{2\sqrt{f}} \\ 
0 & \frac{1}{2\sqrt{h}}
\end{array} 
\right)\V{d}=
\dfrac{1}{4\sqrt{bgfh}}\left(
\begin{array}{c}
\sqrt{2g\sqrt{fh}} \\ 
\sqrt{bh} \\ 
\sqrt{bf}
\end{array} 
\right),\quad
\V{y}^1_{1/2}=(\mbox{need not to calculate}), \quad
\V{y}_0^2=\V{e}_1.
\end{gather*}
Then, the quintuple $\mathcal{X}_+=(Z;a+\mu t^{1/2};\{\V{x}^1,t^\emptyset\V{x}^2\},\{\V{y}^1,t^\emptyset\V{y}^2\})$, where $\V{x}^1:=\V{x}^1_0+\V{x}^1_{1/2}t^{1/2}$, $\V{x}^2:=\V{x}^2_0$,
$\V{y}^1:=\V{y}^1_0+\V{y}^1_{1/2}t^{1/2}$ and $\V{y}^2:=\V{y}^2_0$, is a EQ with Perron-Frobenius property, and the pair $(Z;\mathcal{X}_+)$ is a PF-data.
Because $\mathrm{rk}(\mathcal{X}_+)=1$, we can extend $\mathcal{X}_+$ to any greater depth by the methods in \S \ref{rk1}.
Finally, we have the Perron root $\lambda=a+\mu t^{1/2}+\nu t+\cdots$ ($\nu=\frac{1}{4\sqrt{bf}}(k-\frac{bg}{2h})$) and corresponding Perron vector
\[
\V{x}=
\left(
\begin{array}{c}
\sqrt{2b\sqrt{fh}} \\ 
\sqrt{fg} \\ 
\sqrt{hg}
\end{array} 
\right)+
\left(
\begin{array}{c}
\alpha \\ 
0 \\ 
\beta
\end{array} 
\right)t^{1/2}+\cdots,\qquad
\left(
\alpha=\sqrt{2b\sqrt{fh}}(1-\frac{\nu}{\mu})+\frac{b\sqrt{fg}}{\mu},\quad
\beta=g\sqrt{\frac{b}{2\sqrt{fh}}}
\right).
\]

\subsubsection{The case if $a>\sqrt{fh}$.}
Recall that $G(\mathfrak{A})=\left(
\vcenter{
\xymatrix{
\alpha_1 \ar@<0.5ex>[r]^{1/2}& \alpha_2 \ar@<0.5ex>[l]^{1/2}
}
}
\right)$.
According to the inequality $a>\sqrt{fh}$, we define $E=\{\alpha_1\}$.

Through the similarity transformation $\mathcal{F}:=S_{\alpha_1}(-1/4)\cdot G(\mathfrak{A})=
\left(
\vcenter{
\xymatrix{
\alpha_1 \ar@<0.5ex>[r]^{1/4}& \alpha_2 \ar@<0.5ex>[l]^{3/4}
}
}
\right)
$, we obtain the $E$-gently-slanted digraph $\mathcal{F}$.
Let $F$ be the digraph defined by the diagram $\begin{array}{ccc}
G & \sim & F\\
\downarrow & & \downarrow\\
G(\mathfrak{A}) & \sim & \mathcal{F}
\end{array}$.
(To be more precise, $F:=\left(
\vcenter{
\xymatrix{
1 \ar@(ld,lu)^{0}\ar@<0.5ex>[r]^{9/4}\ar@<0.5ex>[d]^{1/4}& 2 \ar@<0.5ex>[l]^{3/4}\ar@(rd,ru)_3 \ar@<0.5ex>@/^15pt/[dl]^0\\
3 \ar@<0.5ex>[u]^{7/4}\ar@<0.5ex>@/_15pt/[ur]^0\ar@(ld,lu)^1
}
}
\right)
$.)

Set $W:=Y^{F/G}=\left(
\begin{array}{ccc}
a & bt^{3/4} & ct^{7/4}\\
dt^{9/4} & et^3 & f\\
gt^{1/4} & h & kt
\end{array}
\right)$.
Then, $(W;\mathcal{X}(\alpha_1),\mathcal{X}(\alpha_2))$ is a PF-data.
Because $a>\sqrt{fh}$, we have $(a-W[\alpha_2])^{-1}\cdot 
\left(
\begin{array}{c}
dt^{9/4} \\
gt^{1/4}
\end{array}
\right)=\dfrac{1}{a}
\left(
1+\dfrac{W[\alpha_2]}{a}+\dfrac{(W[\alpha_2])^2}{a^2}+\cdots
\right)
\left(
\begin{array}{c}
dt^{9/4} \\
gt^{1/4}
\end{array}
\right)=
\left(
\begin{array}{c}
\frac{gf}{a^2} \\
\frac{g}{a}
\end{array}
\right)t^{1/4}+\cdots
>0$.
Therefore, the quintuple $\mathcal{X}=(W;a;\{\V{x}^1,t^\emptyset\V{x}^2\},\{\V{y}^1,t^\emptyset\V{y}^2\})$, where 
$$
\V{x}^1={}^T(1,\frac{fg}{a^2}t^{1/4},\frac{g}{a}t^{1/4}),\quad \V{x}^2={}^T(0,\sqrt{f},\sqrt{h}),\quad \V{y}^1={}^T(1,0,0)+\cdots,\quad \V{y}^2=\frac{1}{2\sqrt{fh}}{}^T(0,\sqrt{h},\sqrt{f})
$$ 
is a EQ with Perron-Frobenius property, and the pair $(W;\mathcal{X})$ is a PF-data of depth $1/4$.

\subsubsection{The case if $a<\sqrt{fh}$.}
In this case $E=\{\alpha_2\}$.
Let $\mathcal{F}:=S_{\alpha_2}(-1/4)\cdot G(\mathfrak{A})=
\left(
\vcenter{
\xymatrix{
\alpha_1 \ar@<0.5ex>[r]^{3/4}& \alpha_2 \ar@<0.5ex>[l]^{1/4}
}
}
\right)
$.
Define $F$ by the diagram $\begin{array}{ccc}
G & \sim & F\\
\downarrow & & \downarrow\\
G(\mathfrak{A}) & \sim & \mathcal{F}
\end{array}$.
Set $W:=Y^{F/G}=\left(
\begin{array}{ccc}
a & bt^{1/4} & ct^{5/4}\\
dt^{11/4} & et^3 & f\\
gt^{3/4} & h & kt
\end{array}
\right)$.
Then, $(W;\mathcal{X}(\alpha_1),\mathcal{X}(\alpha_2))$ is a PF-data.
Because $a<\sqrt{fh}$, we have $(\sqrt{fh}-a)^{-1}(bt^{1/4},ct^{5/4})\cdot
\left(
\begin{array}{c}
\sqrt{f}\\
\sqrt{h}
\end{array}
\right)=\dfrac{b\sqrt{f}t^{1/4}}{\sqrt{fh}-a}+\cdots>0$.
The quintuple $\mathcal{X}=(W;\sqrt{fh};\{\V{x}^1,t^\emptyset\V{x}^2\},\{\V{y}^1,t^\emptyset\V{y}^2\})$, where $$
\V{x}^1={}^T(\frac{b\sqrt{f}}{\sqrt{fh}-a}t^{1/4},\sqrt{f},\sqrt{h}),\quad \V{x}^2=\V{e}_1,\quad \V{y}^1=\frac{1}{2\sqrt{fh}}{}^T(0,\sqrt{g},\sqrt{f})+\cdots,\quad\V{y}^2=\V{e}_1
$$
is a EQ with Perron-Frobenius property, and the pair $(W;\mathcal{X})$ is a PF-data of depth $1/4$.

\subsection{Example II}\label{sec7.2}

In most cases, our method is applicable for non-negative Puiseux matrices.
Let 
\[
Y=\left(
\begin{array}{ccccc}
\lambda_0 & 1 &  &  &  \\ 
t^{A_1} & \lambda_0 & 1 &  &  \\ 
 & t^{A_2} & \lambda_0 & \ddots &  \\ 
 &  & \ddots & \ddots & 1 \\ 
 &  &  & t^{A_{N-1}} & \lambda_0
\end{array} 
\right),\quad 0<A_1<A_2<\dots<A_{N-1},\quad \lambda_0>0.
\]
Then, the digraph $G:=G(Y)$ is expressed as
$
\left(
\xymatrix{
1\ar@(ld,lu)^0\ar@<0.5ex>[r]^{A_1}& 2 \ar@<0.5ex>[l]^{0}\ar@(ru,lu)_0\ar@<0.5ex>[r]^{A_2} & 3 \ar@<0.5ex>[l]^0\ar@(ru,lu)_0\ar@<0.5ex>[r]^{A_3} &\cdots\ar@<0.5ex>[l]^0 \ar@<0.5ex>[r]^{A_{N-1}}& N\ar@<0.5ex>[l]^0\ar@(rd,ru)_0
}
\right)
$.
Define $H:=\{S_2(\alpha_2)\circ S_3(\alpha_3)\circ\dots\circ S_N(\alpha_{N})\}\cdot G$, where $\alpha_n=\frac{3-n}{2}A_1+(A_2+A_3+\dots+A_{n-1})-(n-2)\epsilon$ and $\epsilon>0$.
The digraph $H$ is expressed as
$
\left(
\xymatrix{
1\ar@(ld,lu)^0\ar@<0.5ex>[r]^{\frac{A_1}{2}}&\ \ \  2\ \ \  \ar@<0.5ex>[l]^{\frac{A_1}{2}}\ar@(ru,lu)_0\ar@<0.5ex>[r]^{\frac{A_1}{2}+\epsilon} & \ \ \ 3\ \ \  \ar@<0.5ex>[l]^{A_2-\frac{A_1}{2}-\epsilon}\ar@(ru,lu)_0\ar@<0.5ex>[r]^{\frac{A_1}{2}+\epsilon} &\ \ \cdots\ \ \ar@<0.5ex>[l]^{A_3-\frac{A_1}{2}-\epsilon} \ar@<0.5ex>[r]^{\frac{A_1}{2}+\epsilon}& \ \ \ N \ar@<0.5ex>[l]^{A_{N-1}-\frac{A_1}{2}-\epsilon}\ar@(rd,ru)_0
}
\right)$.
The digraph $H$ is a flat-slanted whose leading term is a disjoint union of strongly connected components $\alpha_1=\{1,2\}$, $\alpha_2=\{3\},\dots,\alpha_{N-1}=\{N\}$.
Let 
$$
Z:=Y^{H/G}=
\left(
\begin{array}{ccccc}
\lambda_0 & t^\frac{A_1}{2} &  &  &  \\ 
t^{\frac{A_1}{2}} & \lambda_0 & t^{A_2-\frac{A_1}{2}-\epsilon} &  &  \\ 
 & t^{\frac{A_1}{2}+\epsilon} & \lambda_0 & \ddots &  \\ 
 &  & \ddots & \ddots & t^{{A_{N-1}-\frac{A_1}{2}-\epsilon}} \\ 
 &  &  & t^{\frac{A_1}{2}+\epsilon} & \lambda_0
\end{array} 
\right).
$$ 
Consider the quintuples
$$
\mathcal{X}[\alpha_1]=\left(Z[\alpha_1];\lambda_0+t^{\frac{A_1}{2}};
\left\{
\left(
\begin{array}{c}
1 \\ 
1 
\end{array} 
\right),
t^\emptyset
\left(
\begin{array}{c}
0 \\ 
1 
\end{array} 
\right)
\right\},
\left\{
\left(
\begin{array}{c}
\frac{1}{2} \\ 
\frac{1}{2}
\end{array} 
\right),
t^\emptyset
\left(
\begin{array}{c}
0 \\ 
1 
\end{array} 
\right)
\right\}\right)
$$ 
and $\mathcal{X}[\alpha_n]:=(\lambda_0;\lambda_0;\{1\},\{1\})$ ($n\geq 2$).
Then, $\mathfrak{A}:=(Z;\mathcal{X}(\alpha_1),\dots,\mathcal{X}(\alpha_N))$ is a PF-data.
Because $\lambda_0+t^{\frac{A_1}{2}}$ is strictly greater than $\lambda_0$, $\mathfrak{A}$ is non-singular.

Let $E=\{\alpha_1\}$.
The adjacency matrix $G(\mathfrak{A})$ is expressed as
$\left(
\xymatrix{
\alpha_1\ \ \  \ar@<0.5ex>[r]^{\frac{A_1}{2}+\epsilon} & \ \ \ \alpha_2\ \ \  \ar@<0.5ex>[l]^{A_2-\frac{A_1}{2}-\epsilon}\ar@<0.5ex>[r]^{\frac{A_1}{2}+\epsilon} &\ \ \cdots\ \ \ar@<0.5ex>[l]^{A_3-\frac{A_1}{2}-\epsilon} \ar@<0.5ex>[r]^{\frac{A_1}{2}+\epsilon}& \ \ \ \alpha_{N-1} \ar@<0.5ex>[l]^{A_{N-1}-\frac{A_1}{2}-\epsilon}
}
\right)
$.
Note that $G(\mathfrak{A})$ is already  a $E$-gently-slanted.

Let $\beta_1:=\{\alpha_1\}$, $\beta_2:=\{\alpha_2,\dots,\alpha_{N-1}\}$
If $\epsilon$ is sufficiently small, we have 
$$
((\lambda_0+t^{\frac{A_1}{2}})-Z[\beta_2])^{-1}\left(
\begin{array}{cc}
0 & t^{\frac{A_1}{2}+\epsilon} \\ 
\vdots & \vdots\\
0 & 0
\end{array} 
\right)
\left(
\begin{array}{c}
1 \\ 
1
\end{array} \right)\equiv
\left(
\begin{array}{c}
t^{\epsilon} \\ 
t^{2\epsilon} \\ 
\vdots \\ 
t^{(N-2)\epsilon}
\end{array} 
\right)\mod{I((N-2)\epsilon)}.
$$
Let 
\begin{gather*}
\V{x}^1:={}^T(1,1,t^{\epsilon},t^{2\epsilon},\dots,t^{(N-2)\epsilon}),\quad
\V{x}^2:=\V{e}_2t^{(N-2)\epsilon},\dots,\quad
\V{x}^N:=\V{e}_Nt^{(N-2)\epsilon},\\
\V{y}^1:={}^T(\tfrac{1}{2},\tfrac{1}{2},0,\dots,0)+\cdots,\quad 
\V{y}^2:=\V{e}_2t^{(N-2)\epsilon},\dots,\quad
\V{y}^N:=\V{e}_Nt^{(N-2)\epsilon}.
\end{gather*}
By the arguments in \S \ref{rk1}, the quintuple $\mathcal{X}=(Z;\lambda_0+t^{\frac{A_1}{2}};\{\V{x}^1,t^\emptyset\V{x}^i\},\{\V{y}^1,t^\emptyset\V{y}^i\})$ is an EQ with Perron-Frobenius property. 
Finally we obtain the PF-data $(Z;\mathcal{X})$ of depth $(N-2)\epsilon$.

\section{Conclusion and future problems}

In this paper, we have constructed a combinatorial method to calculate the Perron roots/vectors of positive Puiseux matrices.
Our method is applicable for the large class of positive Puiseux matrices with the {\bf condition of genericness} (\S \ref{sec.generic}).
We expect our result to evoke new problems about the structure of the ring of real Puiseux matrices, tropical matrices, Puiseux matrices over a real closed field, {\it etc}.
The following is a list of future problems:

\subsection*{Non-negative real Puiseux matrices}

As introduced in \S \ref{sec7.2}, the Perron root/vector of non-negative matrix can be calculated  by our method in most cases. 
It should be interesting to find an algorithm to calculate them for general non-negative matrices.

\subsection*{Weak Perron-Frobenius property}

A real (might be negative) square matrix is said to possess the {\it weak Perron-Frobenius property} \cite{johnson2004matrices} if its spectral radius is a positive eigenvalue corresponding to non-negative left and right eigenvectors.
Our method can construct many examples of these matrices.
For example, the matrix $Y$ in Example \ref{ex2.16} possesses the weak Perron-Frobenius property for $t=-\epsilon$ ($0<\epsilon\ll 1$).
Many of topological properties of the set of matrices with the weak Perron-Frobenius property have been derived by perturbation methods \cite{elhashash2008general}.
It is expected that more detailed shape of this set will be determined.

\appendix

\section{The Tarski-Seidenberg principle}\label{secaa}

\subsection{General facts about real closed fields}

In this section, several basic facts about real close fields are introduced.
The contents of this section are based on the textbook \cite[\S 1,2,3]{basu2006algorithms}.

A field $k$ is said to be {\it real closed} if (i) there exists a total order $\geq$ so that $k=(k,\geq)$ to be an ordered field, and (ii) the intermediate value theorem holds for all polynomials over $k$.
The real field $\RR$ is a typical example of real closed fields.
Any real closed field contains $\RR$ as a subfield.
An element $x\in k$ is {\it positive} if $x>0$.

Let $k$ be a real closed field.
It is known that $\sqrt{-1}\notin k$. 
Especially, $k$ is not algebraically closed. 
It is also well-known that the field $K:=k+k\sqrt{-1}$ is algebraically closed.
By the intermediate value theorem, any positive element in $k$ is a square element.
Therefore, the sign $\sqrt{\phantom{2}}$ has a meaning over $k$.

Any element of $x\in K$ is expressed as $x=\alpha+\beta\sqrt{-1}$ ($\alpha,\beta\in k$) uniquely.
We define the {\it norm of} $x$ by $\zet{x}:=\sqrt{\alpha^2+\beta^2}$.

Let $D$ be an ordered ring.
We define the {\it language of ordered fields with coefficients in $D$} as follows \cite[p.58]{basu2006algorithms}.
By induction, we give the definition of {\it formulas} and the set $\mathrm{Free}(\Phi)$ of {\it free variables of formula} $\Phi$:
\begin{enumerate}
\item An {\it atom} is ($P>0$) or ($P=0$), where $P\in D[X_1,\dots,X_k]$.
An atom is a formula. 
Set $\mathrm{Free}(P>0)=\mathrm{Free}(P=0):=\{X_1,\dots,X_k\}$.
\item If $\Phi_1$ and $\Phi_2$ are formulas, then $\Phi_1\wedge \Phi_2$ and $\Phi_1\vee \Phi_2$ are formulas.
Define $\mathrm{Free}(\Phi_1\wedge \Phi_2)=\mathrm{Free}(\Phi_1\vee \Phi_2):=\mathrm{Free}(\Phi_1)\cup\mathrm{Free}(\Phi_2)$.
\item If $\Phi$ is a formula, then $\neg\Phi$ is a formula with $\mathrm{Free}(\neg \Phi):=\mathrm{Free}(\Phi)$.
\item If $\Phi$ is a formula and $X\in\mathrm{Free}(\Phi)$, then $(\forall X)\Phi$ and $(\exists X)\Phi$ are formulas with $\mathrm{Free}((\forall X) \Phi)=\mathrm{Free}((\exists X)\Phi):=\mathrm{Free}(\Phi)\setminus\{X\}$.
\end{enumerate}
For two formulas $\Phi_1$ and $\Phi_2$, $\Phi_1\to \Phi_2$ is the formula $(\neg \Phi_1)\vee \Phi_2$.
A formula is said to be {\it quantifier free} if it contains no quantifier ($\forall$, $\exists$).
A {\it sentence} is a formula without free variable.

Let $\Phi$ be a formula in the language of ordered fields with coefficients in $D$ and $k$ be a ordered field such that $D\subset k$.
The {\it realization of $\Phi$ in $k$} is the formula which is obtained from $\Phi$ by replacing "$\forall X$" to "$\forall X\in k$" and "$\exists X$" to "$\exists X\in k$".
Denote $\Phi\vert_k$ the realization of $\Phi$ in $k$.
A sentence $\Phi$ is said to be {\it true in $k$} if $\Phi\vert_k$ is true.
A {\it proof of $\Phi_n$ in $k$} is a sequence of formulas $(\Phi_1,\dots,\Phi_n)$ such that, for all $i$, (i) $\Phi_i\vert_k$ is an axiom of real closed field $k$ or (ii) $\Phi_i\vert_k$ can be derived from $\Phi_1\vert_k,\dots,\Phi_{i-1}\vert_k$ by inference rule.
A formula $\Phi$ is {\it provable in $k$} if there exists a proof of $\Phi$ in $k$.

The following theorems are quite important \cite[Theorems 2.77 and 2.80]{basu2006algorithms}:
\begin{thm}[Quantifier elimination]\label{quanti}
Let $k$ be a real closed field and $D$ be an ordered subring of $k$.
If a sentence $\Phi$ in the language of ordered fields with coefficients in $D$ is true in $k$, then it is provable in $k$.
\end{thm}

\begin{thm}[Tarski-Seienberg principle]\label{tarski}
Let $k$ and $k'$ be real closed fields and $D$ be an ordered ring such that $D\subset k\subset k'$.
If $\Phi$ is a sentence in the language of ordered fields with coefficients in $D$, then
\begin{center}
$\Phi$ is true in $k$ $\iff$ $\Phi$ is true in $k'$.
\end{center}
\end{thm}
\begin{cor}
Let $\Phi$ be a sentence in the language of ordered fields with coefficients in $\ZZ$.
If $\Phi$ is true in $\RR$, then $\Phi$ is true in any real closed field.
\end{cor}

\subsection{An example how to use the Tarski-Seidenberg principle}

The following is a typical example of the application of the Tarski-Seidenberg principle.
\begin{lemma}\label{lemmaa4}
Let $X\in k^{N\times N}$ be a square matrix all the entries of which are nonnegative.
Denote $\lambda_1,\dots,\lambda_N\in K$ the eigenvalues of $X$.
If $\mu>\zet{\lambda_n}$ for all $n$, all the entries of the matrix $(\mu\cdot \mathrm{Id}-X)^{-1}$ are nonnegative.
\end{lemma}

\proof Let $X_{i,j}$, $Y_{i,j}$ ($i,j=1,\dots,N$), $A$ and $B$ denote intermediates.
Let $f(A,B,X),g(A,B,X)\in \ZZ[A,B,X_{i,j}]$ be polynomials such that $\det((A+B\sqrt{-1})\mathrm{Id}-X)=f(A,B,X)+g(A,B,X)\sqrt{-1}$.
Consider the formula $\Phi$ which is described as follows:
\begin{align*}
&\forall X_{i,j}, \forall Y_{i,j},\forall \mu, \forall A, \forall B\\
&
\left(
\bigwedge_{i,j}(X_{i,j}\geq 0)\wedge
\left(
(
f(A,B,X)=0
\wedge 
g(A,B,X)=0
)\to (\mu^2 >A^2+B^2)
\right)\wedge
\bigwedge_{i,k}
(
\sum_{j=1}^{N}Y_{i,j}(\mu\delta_{j,k}-X_{j,k})=\delta_{i,k}
)
\right)\\
&\to \bigwedge_{i,j}(Y_{i,j}\geq 0).
\end{align*}
The realization of $\Phi$ in $k$ is nothing but the statement of the lemma.
Because $\Phi$ is a sentence in the language of ordered fields with coefficients in $\ZZ$, $\Phi$ is true in $k$ if and only if $\Phi$ is true in $\RR$.
As it is soon proved that $\Phi$ is true in $\RR$, $\Phi$ is also true in any $k$.  $\qed$

\section*{Acknowledges} 
The author is grateful to Kyo Nishiyama for his insightful comments on this manuscript.
This work was supported by KAKENHI (26800062).

\bibliography{test}

\end{document}